\documentclass[11pt]{amsart}
\usepackage{amsmath}
\usepackage{amsthm,amssymb,amsfonts,dsfont,leftindex} 
\usepackage{geometry}
\usepackage{tikz-cd} 							
\usepackage{varioref}
\usepackage{hyperref}
\hypersetup{colorlinks,linkcolor=magenta,citecolor=teal,linktocpage}
\usepackage[noabbrev, capitalise, nameinlink]{cleveref}
\usepackage[alphabetic]{amsrefs} 					
\usepackage{thmtools}								
\usepackage{style}
\begin{document}	
\setcounter{tocdepth}{1}
\title[A Beauville-Laszlo-type descent theorem]{A Beauville-Laszlo-type descent theorem for locally Noetherian schemes}
\author{Robin Louis}

\begin{abstract}
Let $X$ be a locally Noetherian scheme with a closed subscheme $Z$. Let $\X$ be the completion of $X$ at $Z$, considered as a formal scheme.
We show that a coherent sheaf on $X$ is equivalently given by a coherent sheaf on $\X$, a coherent sheaf on the complement of $Z$, and an isomorphism of pullbacks of these sheaves to a certain adic space $W$.
By defining $W$ as an adic space instead of as a Berkovich space we are able to generalize the descent result of \cite{BT} from finite type $k$-schemes to locally Noetherian schemes.
\end{abstract}
\maketitle

\tableofcontents
\newpage
\section*{Introduction}
The present paper deals with a descent result which is motivated by the Beauville-Laszlo theorem. Therefore it is illuminating to first recall this theorem to understand where the current paper is headed.

Let $A$ be a ring and let $f\in A$ be an element. Let $\comp{A}$ denote the $f$-adic completion of $A$.
For an $A$-algebra $B$ let $\Modf[B]$ denote the full subcategory of $\Mod[B]$ of $f$-regular modules.
In \cite{BL} Beauville and Laszlo showed that if $f$ is a nonzerodivisor the following square given by basechange is cartesian:
\[
\begin{tikzcd}
	\Modf[A]\arrow[r]\arrow[d]&\Mod[{A[f^{-1}]}]\arrow[d]\\
	\Modf[\comp{A}]\arrow[r]&\Mod[\comp{A}{[f^{-1}]}]
\end{tikzcd}
\]
Denoting $X=\Spec{A}$, $Z=V(f)$, $U=X$\,\textbackslash\,$Z$ and the completion of $X$ at $Z$ by $\X$ we can alternatively write this square as
\[
\begin{tikzcd}
	\QCohf[X]\arrow[r,"j^*"]\arrow[d,"i^*",swap]&\QCoh[U]\arrow[d]\\
	\QCohf[\X]\arrow[r]&\Mod[\comp{A}{[f^{-1}]}]\mathrlap{,}
\end{tikzcd}
\]
where the left and top functor are pullback along $i\colon\X\to X$ and $j\colon U\to X$.
When trying to generalize this to a global situation for a general scheme $X$ with some closed subscheme $Z$ it is thus clear that some objects over $X$, $U$ and $\X$ should be involved (there is the additional technical problem that one needs to be careful with quasicoherent sheaves on non-affine formal schemes and with non-Noetherian formal schemes in general, see the discussion in \cite{duality}); however it is not as clear what a good generalization of $\comp{A}[f^{-1}]$ would be.

This is a form of descent statement, so a first idea might be to consider the pullback $\X\times_XU$ as formal schemes, but $\X$ has underlying topological space $Z$, so this pullback is empty.

However, if one is satisfied with a statement for coherent sheaves and assuming that $X$ is locally Noetherian there exists a good replacement for $\X$ and $U$ which we can naively take the pullback of. These assumptions resolve the technical problem as well, as coherent sheaves on locally Noetherian formal schemes are well-behaved.
Using the setup of adic spaces as introduced by Huber one can consider the adification $X^{\ad}$ of $X$ as well as the adification $\X^{\ad}$ of $\X$.
As $X$ is a scheme we have a morphism of locally ringed spaces
\[
s_X\colon (X^{\ad},\mathcal{O}_{X^{\ad}})\to(X,\mathcal{O}_X),
\]
which locally takes a valuation to its support.
Let $\U$ be the open adic subspace of $X^{\ad}$ on the preimage of $U$ under this map.
When replacing $X$ with $X^{\ad}$, $\X$ with $\X^{\ad}$ and $U$ with $\U$ the wanted geometric replacement of $\comp{A}[f^{-1}]$ is just $W\coloneqq \X^{\ad}\times_{X^{\ad}}\U$.

It turns out that in the locally Noetherian setup one does not need any assumption on the ideal cutting out $Z$ or the modules, so the theorem that we will prove is the following:

\begin{theorem*}\label{mainthm}
	There exist functors $\Coh[U]\to \Coh[W]$ and $\Coh[\X]\to\Coh[W]$ making the square
	\[
	\begin{tikzcd}
		\Coh[X]\arrow[r,"j^*"]\arrow[d,"i^*",swap]&\Coh[U]\arrow[d]\\
		\Coh[\X]\arrow[r]&\Coh[W]
	\end{tikzcd}
	\]
	into a cartesian square.
	These functors are given by pulling back coherent sheaves along a morphism of the underlying locally ringed spaces.
\end{theorem*}
We will also show that the square
\[
\begin{tikzcd}
	\Coh[X^{\ad}]\arrow[d]\arrow[r]&\Coh[\U]\arrow[d]\\
	\Coh[\X^{\ad}]\arrow[r]&\Coh[W]
\end{tikzcd}
\]
given by pulling back along the morphisms defining $W$ is cartesian.

In the affine principal ideal case this statement recovers a version of the Beauville-Laszlo theorem, since the above square is then given by
\[
\begin{tikzcd}
	\Coh[\Spa{A}{A}]\arrow[r]\arrow[d]&\Coh[\Spa{A[f^{-1}]}{A^{\cl}}]\arrow[d]\\
	\Coh[\Spa{A_f}{A}]\arrow[r]&\Coh[\Spa{A_f[f^{-1}]}{A^{\cl}}]\mathrlap{,}
\end{tikzcd}
\]
which is equivalent to
\[
\begin{tikzcd}
	\fMod[A]\arrow[r]\arrow[d]&\fMod[A{[f^{-1}]}]\arrow[d]\\
	\fMod[\comp{A}]\arrow[r]&\fMod[\comp{A}{[f^{-1}]}]
\end{tikzcd}
\]
by taking global sections.
One might think that this last square being cartesian follows immediately from faithfully flat descent since $\comp{A}$ is flat over $A$ if $A$ is Noetherian.
This is not the case however, since faithfully flat descent involves a descent datum over $\comp{A}\otimes_A\comp{A}$ and one over $\comp{A}\otimes_AA[f^{-1}]$, i.e.\@ $\comp{A}[f^{-1}]$, while the statement above only has an isomorphism over the latter.
That the isomorphism over $\comp{A}\otimes_A\comp{A}$ is not needed to construct the module over $A$ is not clear a priori.

The idea of using a rigid analog of $\X$ and $U$ to show a Beauville-Laszlo-style theorem for non-affine schemes is not new, for example in \cite{BT} Ben-Bassat and Temkin have already shown a variant of the present generalization in a setup using Berkovich spaces when assuming that $X$ is a $k$-scheme of finite type.
They also suggest a generalization of their result using adic spaces in Part 4.6 of their paper, however the present proof proceeds somewhat differently than the proof in \cite{BT}.

The proof of the theorem will be accomplished in several steps after we recall some basics about coherent sheaves on adic spaces in the first section.
In the second section we will define the functors occurring in the theorem and show that the square is commutative, which enables us to properly formulate the main theorem as Theorem~\ref{mainthmprop}.
Next we will prove the Key~Lemma~\ref{keylemma}, which is both the main step towards the proof of the theorem in the case where $Z$ is locally principal, as well as enabling us to reduce to this case by blowing up $Z$.
In the fourth section we compile the steps taken so far to first prove the theorem in the locally principal case and then in the general case.
The last section deals with a discussion of further corollaries and the explicit description of $W$.
\subsection*{Notation and Conventions}
If $A$ is a ring and $I$ is an ideal of $A$, then $A$ considered as a topological ring with the $I$-adic topology is denoted by $A_I$.
If the ideal that topologizes $A$ is clear from the context we will drop the $I$ for brevity.

Let us now quickly talk about fibre products and commutative squares in categories.
Consider a square of functors
\[
\begin{tikzcd}
	\catC\arrow[r,"G^\prime"]\arrow[d,"F^\prime",swap]&\catD\arrow[d,"F"]\\
	\catE\arrow[r,"G"]&\catF\mathrlap{.}
\end{tikzcd}
\]
Let $\sigma$ be a natural isomorphism between the two compositions $\catC\to\catF$.
We say that the commutativity of the square up to natural isomorphism is exhibited by $\sigma$.

We denote the fibre product (one might more accurately call this a strict iso-comma object in $\catname{CAT}_{2,1}$) of a cospan of categories as above by $\catD\times_{\catF}\catE$.
An object $(b,c,\varphi)$ of this category is an object $b$ of $\catD$ and an object $c$ of $\catE$ together with an isomorphism $\varphi\colon F(b)\to G(c)$ in $\catF$.
A morphism $(b,c,\varphi)\to(b^\prime,c^\prime,\varphi^\prime)$ is given by a morphism $f\colon b\to b^\prime$ in $\catD$ and a morphism $g\colon c\to c^\prime$ in $\catE$ such that
\[
\begin{tikzcd}
	F(b)\arrow[r,"F(f)"]\arrow[d,"\varphi",swap]&F (b^\prime)\arrow[d,"\varphi^\prime"]\\
	G(c)\arrow[r,"G(g)"]&G (c^\prime)
\end{tikzcd}
\]
commutes.
Forgetting about the isomorphism $\varphi$ and one of $b$ or $c$ yields projection functors
\[
\pr_1\colon\catD\times_{\catF}\catE\to\catD
\]
and
\[
\pr_2\colon\catD\times_{\catF}\catE\to\catE\mathrlap{.}
\]
The isomorphisms given in the data of the objects assemble into a natural isomorphism
\[
\varphi\colon F\circ\pr_1\to G\circ\pr_2\mathrlap{,}
\]
which fulfills the following universal property:

For every category $\catC$ together with functors $G^\prime\colon\catC\to\catD$ and $F^\prime\colon\catC\to\catE$ and a natural isomorphism $\sigma$ of the compositions $\catC\to\catF$ there exists a unique functor $H\colon\catC\to\catD\times_{\catF}\catE$ which makes the triangles in the following diagram
\[
\begin{tikzcd}
	\catC\arrow[ddr,swap,"F^\prime",bend right=30]\arrow[rrd,"G^\prime",bend left=20]\arrow[dr,"H"]\\[-4pt]
	&[-18pt]\catD\times_{\catF}\catE\arrow[r]\arrow[d]&\catD\arrow[d,"F"]\\
	&\catE\arrow[r,"G"]&\catF\mathrlap{.}
\end{tikzcd}
\]
commute (not up to natural isomorphism, but honestly commute) and such that the induced natural transformation $\varphi\circ H$ between the compositions $\catC\to\catF$ is equal to $\sigma$.

Then we say that a square whose commutativity is exhibited by $\sigma$ is cartesian if the unique induced functor
\[
\catC\to\catD\times_{\catF}\catE
\]
is an equivalence.

We will in the following text often abusively call a square of functors commutative if there exists some canonical natural isomorphism exhibiting the commutativity.
This isomorphism will however usually be clear from the context or be provided by lemmata that are applied to obtain the commutative square in the first place.
Similarly we will call a square that commutes up to a canonical natural isomorphism cartesian if it is cartesian with regard to this natural isomorphism.

\subsection*{Acknowledgements and Thanks}
I am very grateful to Manuel Hoff and Mingjia Zhang for their help in mathematical, proofreading and typesetting questions.

This paper is a slightly adapted version of my Master's thesis at the University of Bonn.
Thus I am particularly indebted to my advisor Peter Scholze for his support and the informative and enlightening discussions I had with him during the development of the thesis.
Moreover I want to thank Johannes Ansch\"utz for being willing to co-correct the thesis.

Special thanks are given to my family for their unwavering aid, for their encouragement and for everything else.
\newpage
\section{Generalities}
The precise formulation of the main theorem will require basic facts about the theory of adic spaces and coherent sheaves on adic spaces. Most of these can be found in \cite{Hu} and we will assume a certain familiarity with adic spaces themselves, but we will give a short recap of the most relevant points regarding adifications, fibre products, coherent sheaves and pullbacks here.
\subsection{Adification}
\begin{lemma}
	Let $\X$ be a locally Noetherian formal scheme. There exists an adic space $\X^{\ad}$ called the \emph{adification of} $\X$ and a morphism of topologically locally ringed spaces
	\[
	p_{\X}\colon (\X^{\ad},\mathcal{O}_{\X^{\ad}}^+)\to (\X,\mathcal{O}_{\X})
	\]
	called the \emph{structure morphism}, such that for every adic space $T$ and every morphism of topologically locally ringed spaces
	\[
	f\colon(T,\mathcal{O}_T^+)\to(\X,\mathcal{O}_{\X})
	\]
	there exists a unique morphism $T\to\X^{\ad}$ of adic spaces such that the induced morphism of topologically locally ringed spaces
	\[
	(T,\mathcal{O}_T^+)\to(\X^{\ad},\mathcal{O}_{\X^{\ad}}^+)
	\]
	factors $f$ through $p_{\X}$.
\end{lemma}
\begin{proof}
	This is \cite[Proposition 3.9.13]{Hu}.
\end{proof}
\begin{remark}\label{pglobalsec}
	\item
	\begin{itemize}
		\item 
		The adification assembles into a functor $(-)^{\ad}$ from the category of locally Noetherian formal schemes to the category of adic spaces, where the action on morphisms sends a map of formal schemes
		\[
		f\colon(\X,\mathcal{O}_\X)\to(\Y,\mathcal{O}_\Y)
		\]
		to the morphism
		\[
		f^{\ad}\colon\X^{\ad}\to\Y^{\ad}
		\]
		of adic spaces induced by the universal property, i.e.\@ the unique morphism between the adifications that makes the following square commutative:
		\[
		\begin{tikzcd}
			(\X^{\ad},\mathcal{O}_{\X^{\ad}})\arrow[r,"f^{\ad}"]\arrow[d,"p",swap]&(\Y^{\ad},\mathcal{O}_{\Y^{\ad}})\arrow[d,"p"]\\
			(\X,\mathcal{O}_\X)\arrow[r,"f"]&(\Y,\mathcal{O}_\Y)
		\end{tikzcd}
		\]
		\item 
		The adification functor sends $\Spf{A}$ to $\Spa{A}{A}$. A morphism on formal spectra that is induced by a  continuous morphism $\varphi\colon A\to B$ is sent to the morphism on $\operatorname{Spa}$ which is induced by $\varphi$.
		\item 
		Huber shows in \cite[Lemma 3.9.12]{Hu} where this morphism comes from in the affine case.
		Let $\X=\Spf{A}$ be an locally Noetherian affine formal scheme and let $T$ be an adic space.
		Denoting the category of topologically locally ringed spaces by $\TLRS$ the map 
		\begin{align*}
			\quad\quad\quad\TLRS((T,\mathcal{O}^+_T),(\Spf{A},\mathcal{O}_{\Spf{A}}))\to \TopRing(\mathcal{O}_{\Spf{A}}(\Spf{A}),\mathcal{O}^+_T(T))
		\end{align*}
		sending a morphism to the induced map on global sections is an isomorphism.
		When plugging in $T=\X^{\ad}$, which is then given by $\Spa{A}{A}$, the unique preimage of the identity is $p_{\X}$.
		\item 
		For a Noetherian affine formal scheme $\Spf{A}$ the structure morphism sends a continuous valuation $\nu\in\Spa{A}{A}$ to the ideal $\{\nu<1\}$, which is an open prime ideal in $A$.
		\item 
		In the construction of $\X^{\ad}$ one takes the local adifications $\Spa{A}{A}$ and then glues them together.
		The local Noetherianness is needed to make these sheafy. As we will only consider locally Noetherian schemes and their completions, which are again locally Noetherian, we will be able to take adifications of all schemes and formal schemes that appear.
	\end{itemize}
\end{remark}
In general fibre products and coherent sheaves on adic spaces are not that well-behaved. The adic spaces we are interested in (which are all open subspaces of adifications of locally Noetherian formal schemes) however all fulfill a certain local Noetherianity assumption, which makes fibre products and coherent sheaves work nicely, f.e.\@ consider Lemma~\ref{fibreschemeadic} and Lemma~\ref{cohrspace}.
\begin{definition}
	We call an adic space $X$ \emph{locally Noetherian} if the Huber ring $\mathcal{O}_X(V)$ has a Noetherian ring of definition over which it is finitely generated for every open affinoid subspace $V\subseteq X$.
	
	We call a Huber pair $(A,A^+)$ \emph{Noetherian} if $A$ has a Noetherian ring of definition over which it is finitely generated.
\end{definition}
\begin{lemma}\label{Noetherian}
	Let $\X$ be a locally Noetherian formal scheme. Then $\X^{\ad}$ is locally Noetherian.
\end{lemma}
\begin{proof}
	This is part of \cite[Proposition 4.2.ii]{Hu2}.
\end{proof}
\begin{remark}
	It also follows from \cite[Proposition 4.2.ii]{Hu2} that an affinoid adic space $\Spa{A}{A^+}$ is locally Noetherian iff $(A,A^+)$ is Noetherian, justifying the terminology.
\end{remark}
\subsection{Fibre products}
Since the theorem we are interested in concerns a descent situation we need to discuss fibre products of adic spaces.
In general these need not exist (\cite[Chapter 3.10]{Hu} gives assumptions on the morphisms which secure that the fibre product exists). However we can always form the fibre product along an open immersion:
\begin{lemma}\label{pullback}
	Let $f\colon X\to Y$ be a morphism between adic spaces (schemes, formal schemes) and let $U\subseteq Y$ be an open subspace of $Y$.
	Then the inclusion and the restriction of $f$ induce a cartesian square
	\[
	\begin{tikzcd}
		f^{-1}(U)\arrow[r]\arrow[d]&U\arrow[d]\\
		X\arrow[r,"f"]&Y\mathrlap{.}
	\end{tikzcd}
	\]
\end{lemma}
\begin{proof}
	We show the lemma for $X$ and $Y$ adic spaces, the other two are similar.
	We need to check that the square above is cartesian.
	
	Let $T$ be an adic space together with morphisms $g\colon T\to X$ and $h\colon T\to U$ that become equal after mapping to $Y$. In particular the image of $g$ lies in $f^{-1}(U)$.
	The equalities
	\[
	\mathcal{O}_{f^{-1}(U)}=\mathcal{O}_X|_{f^{-1}(U)}
	\]
	as well as
	\[
	\nu_{f^{-1}(U),x}=\nu_{X,x},\text{ for all } x\in f^{-1}(U)
	\]
	hold by definition. Thus the morphism $g$ factors through $f^{-1}(U)$ as a morphism of adic spaces. Because $f^{-1}(U)\to X$ is mono it even factors uniquely through $f^{-1}(U)$. To check that the composition to $f^{-1}(U)$ is $h$, compose with $U\to Y$ and use that this map is mono.
	Thus $f^{-1}(U)$ fulfills the universal property of the pullback.
\end{proof}
One can also take a kind of ``fibre product'' of an adic space with a scheme over a scheme, if some finiteness and Noetherianity assumptions are met. We will abusively also call this object a fibre product or pullback of the two, even though they do not live in the same ambient category.
The following lemma discussing this is taken from \cite[Proposition 3.8]{Hu2}.
\begin{lemma}\label{fibreschemeadic}
	Let $f\colon X^\prime\to X$ be a morphism of locally Noetherian schemes that is locally of finite type.
	Let $S$ be an adic space such that $\mathcal{O}_S(U)$ has a Noetherian ring of definition for every affinoid $U\subseteq S$.
	Let $g\colon (S,\mathcal{O}_S)\to (X,\mathcal{O}_X)$ be a morphism from the underlying locally ringed space of $S$ to $X$.
	
	Then there exists an adic space $S\times_XX^\prime$ together with a morphism $q\colon S\times_XX^\prime\to S$ of adic spaces and a morphism $r\colon (S\times_XX^\prime,\mathcal{O}_{S\times_XX^\prime})\to (X^\prime,\mathcal{O}_{X^\prime})$ of locally ringed spaces making
	\[
	\begin{tikzcd}
		S\times_XX^\prime\arrow[r,"r"]\arrow[d,"q",swap]&X^\prime\arrow[d,"f"]\\
		S\arrow[r,"g"]&X
	\end{tikzcd}
	\]
	commutative as a diagram in the category of locally ringed spaces, such that for every adic space $T$ with morphisms $T\to S$ of adic spaces and $(T,\mathcal{O}_T)\to (X^\prime,\mathcal{O}_{X^\prime})$ of locally ringed spaces there exists a unique morphism $T\to S\times_XX^\prime$ of adic spaces such that
	\[
	\begin{tikzcd}
		T\arrow[dr]\arrow[rrd,bend left=17]\arrow[ddr,bend right=30]\\[-4pt]
		&[-18pt]S\times_XX^\prime\arrow[r,"r"]\arrow[d,"q",swap]&X^\prime\arrow[d,"f"]\\
		&S\arrow[r,"g"]&X
	\end{tikzcd}
	\]
	commutes as a diagram of locally ringed spaces.
	This alternatively means that postcomposition with $q$ and $r$ induces an isomorphism
	\[
	\Ad(T,S\times_XX^\prime)\xrightarrow{\cong}\Ad(T,S)\times_{\TLRS((T,\mathcal{O}_T),X)}\TLRS((T,\mathcal{O}_T),X^\prime)
	\]
	when denoting the category of topologically locally ringed spaces by $\TLRS$.\qed
\end{lemma}
\begin{remark}
	As remarked above all adic spaces that we work with are open subspaces of adifications of locally Noetherian formal schemes. Thus they are all locally Noetherian by Lemma~\ref{Noetherian}, in particular the ring of functions of every affinoid subspace has a Noetherian ring of definition, so we can always form this pullback.
\end{remark}
\begin{remark}\label{openfakepb}
	If $V\to X$ is the inclusion of an open subscheme the fibre product $S\times_XV$ identifies with the open subspace of $S$ on the preimage $g^{-1}(V)$ by considering the universal properties of the two spaces.
\end{remark}
Now we can ask ourselves how one can explicitly describe the fibre product of the map $p$ with an inclusion of an open subscheme $V\hookrightarrow X$. The natural thing to suggest would be that this is just the adification of $V$.
We show a slightly more general statement about a further fibre product with the adification of a completion $\X^{\ad}$.
\begin{lemma}\label{structurepb}
	Let $\X$ be a completion of a locally Noetherian scheme $X$ at a closed subscheme $Z$ and let $V$ be an open subscheme of $X$.
	Denote the fibre product $\X \times_XV$ by $\V$.
	Then we have an isomorphism
	\[
	\V^{\ad}\cong\X^{\ad}\times_{X^{\ad}} (X^{\ad}\times_XV)\mathrlap{.}
	\]
\end{lemma}
\begin{proof}
	Both $\V^{\ad}$ and $\X^{\ad}\times_{X^{\ad}} (X^{\ad}\times_X\V)$ are open subspaces of $\X^{\ad}$ by Lemma~\ref{pullback}.
	The structure map $\V^{\ad}\to\V$ composed with $\V\to V$ and the adification of $\V\to X$ induce a map
	\[
	\V^{\ad}\to X^{\ad}\times_XV\mathrlap{.}
	\]
	This map and the inclusion $\V^{\ad}\to\X^{\ad}$ induce a map
	\[
	\V^{\ad}\to\X^{\ad}\times_{X^{\ad}} (X^{\ad}\times_X\V)\mathrlap{.}
	\]
	We obtain an isomorphism
	\[
	\begin{aligned}
		\Ad(T,\V^{\ad})&\cong\TLRS((T,\mathcal{O}_T^+),\V)\\
		&=\TLRS((T,\mathcal{O}_T^+),\X\times_XV)\\
		&\cong\TLRS((T,\mathcal{O}_T^+),\X)\times_{\TLRS((T,\mathcal{O}_T^+),X)}\TLRS((T,\mathcal{O}_T^+),V)\\
		&\cong\Ad(T,\X^{\ad})\times_{\Ad(T,X^{\ad})}\TLRS((T,\mathcal{O}_T^+),V)\\
		&\cong\Ad(T,\X^{\ad})\times_{\Ad(T,X^{\ad})}\TLRS((T,\mathcal{O}_T^+),X)\times_{\TLRS((T,\mathcal{O}_T),X)}\TLRS((T,\mathcal{O}_T),V)\\
		&\cong\Ad(T,\X^{\ad})\times_{\Ad(T,X^{\ad})}\Ad(T,X^{\ad})\times_{\TLRS((T,\mathcal{O}_T),X)}\TLRS((T,\mathcal{O}_T),V)\\
		&\cong\Ad(T,\X^{\ad})\times_{\Ad(T,X^{\ad})}\Ad(T,X^{\ad}\times_XV)\\
		&\cong\Ad(T,\X^{\ad}\times_{X^{\ad}} (X^{\ad}\times_XV))
	\end{aligned}
	\]
	given by the universal properties of the adification and the fibre product. Going through the isomorphism chain the isomorphism is given by postcomposing with the above map
	\[
	\V^{\ad}\to\X^{\ad}\times_{X^{\ad}}(X^{\ad}\times_X\V)\mathrlap{.}
	\]
	Thus this map is an isomorphism by Yoneda.
\end{proof}
\begin{corollary}\label{schemestructurepb}
	Let $X$ be a locally Noetherian scheme, let $V$ be an open subscheme of $X$.
	Then
	\[
	V^{\ad}\cong X^{\ad}\times_XV\mathrlap{.}
	\]
\end{corollary}
\begin{proof}
	This is Corollary~\ref{structurepb} for $Z=X$.
\end{proof}

\begin{corollary}\label{structurepreimage}
	In the setup of Lemma~\ref{structurepb} the preimage of $\V$ in $\X^{\ad}$ under $p$ is given by $\V^{\ad}$.
\end{corollary}
\begin{proof}
	The square
	\[
	\begin{tikzcd}
		\X^{\ad}\arrow[r]\arrow[d,"p",swap]&X^{\ad}\arrow[d,"p"]\\
		\X\arrow[r]&X
	\end{tikzcd}
	\]
	commutes. The subspace on the preimage of $V$ in $\X^{\ad}$ along the upper composition is
	\[
	\X^{\ad}\times_{X^{\ad}} (X^{\ad}\times_XV)
	\]
	by combining Remark~\ref{openfakepb} and Lemma~\ref{pullback}. By Lemma~\ref{structurepb} this is $\V^{\ad}$.
	On the other hand this is the preimage of $V$ under the lower composition, i.e.\@ the preimage of $\V$ under $p$.
\end{proof}
Now we can show that the intersection of two adifications of open subschemes is the adification of their intersection.
\begin{lemma}\label{openint}
	Let $X$ be a locally Noetherian scheme and let $U$ and $V$ be two open subschemes. Then the adifications of the inclusions induce an isomorphism
	\[
	(U\times_XV)^{\ad}\to U^{\ad}\times_{X^{\ad}}V^{\ad}\mathrlap{.}
	\]
\end{lemma}
\begin{proof}
	For an open subscheme $X^\prime\subseteq X$ the inclusion and structure map induce an isomorphism $X^{\prime\ad}\cong X^{\ad}\times_XX^\prime$ by Corollary~\ref{schemestructurepb}.
	Moreover for $U$ and $V$ open subschemes we have a morphism
	\[
	X^{\ad}\times_X(U\times_XV)\cong X^{\ad}\times_XU\times_{X^{\ad}}X^{\ad}\times_XV
	\]
	given by projecting. This is an isomorphism, since one can check that both sides fulfill the same universal property.
	Thus we have an isomorphism
	\begin{align*}
		(U\times_XV)^{\ad}&\cong X^{\ad}\times_X(U\times_XV)\\
		&\cong X^{\ad}\times_XU\times_{X^{\ad}}X^{\ad}\times_XV\\
		&\cong U^{\ad}\times_{X^{\ad}}V^{\ad}\mathrlap{.}\qedhere
	\end{align*}
\end{proof}
If $X$ is a scheme considered as a formal scheme we have the structure map
\[
p\colon (X^{\ad},\mathcal{O}_{X^{\ad}}^+)\to X\mathrlap{;}
\]
we can also define a different map from the adification to $X$ in the following way:

Let $\Spec{A}$ be an affine scheme.
The identity on global sections induces a morphism of locally ringed spaces
\[
s_{\Spec A}\colon (\Spa{A}{A},\mathcal{O}_{\Spa{A}{A}})\to (\Spec{A},\mathcal{O}_{\Spec{A}})\mathrlap{,}
\]
which sends a valuation $\nu$ on $A$ to its support.
\begin{lemma}
	For a locally Noetherian scheme $X$ the maps $s_{W}$ for $W\subseteq X$ affine open glue to a morphism of locally ringed spaces
	\[
	s_X\colon (X^{\ad},\mathcal{O}_{X^{\ad}})\to (X,\mathcal{O}_X)
	\]
	called the \emph{support morphism}.
\end{lemma}
\begin{proof}
	We can glue the $s_W$ once we check that for affine opens $U$ and $V$ the maps $s_U$ and $s_V$ agree on the intersection $U^{\ad}\cap\, V^{\ad}$ in $X^{\ad}$.
	This intersection is the pullback $U^{\ad}\times_{X^{\ad}}V^{\ad}$ by Lemma~\ref{pullback}, which is $(U\times_XV)^{\ad}$ by Lemma~\ref{openint}.
	Cover $U\times_XV$ with affine opens $\{W_i\}_i$.
	Then $(U\times_XV)^{\ad}$ is covered by the $\{W_i^{\ad}\}_i$. Let $\iota_i\colon W_i\to U$ be the inclusion.
	The morphisms of locally ringed spaces $s_U\circ \iota_i$ and $\iota_i^{\ad}\circ s_{W_i}$ induce the same map on global sections and map to an affine scheme, so they are equal.
	The same argument works for $V$, showing that the restrictions of $s_U$ and $s_V$ to $(U\times_XV)^{\ad}$ are equal, so we can glue the $s_W$ to a morphism of locally ringed spaces
	\[
	s_X\colon (X^{\ad},\mathcal{O}_{X^{\ad}})\to (X,\mathcal{O}_X)\mathrlap{.}\qedhere
	\]
\end{proof}
\begin{remark}
	Later we will with a bit of abusive notation denote $s$ as $s\colon X^{\ad}\to X$ and $p$ as $p\colon \X^{\ad}\to\X$ when considering commutative diagrams with both adic spaces as well as formal schemes.
	These commutative diagrams live in the category of topologically locally ringed spaces, and the forgetful functor from adic spaces has to be chosen accordingly.
	If the diagram contains $p$ the implicit forgetful functor sends an adic space $X$ to $(X,\mathcal{O}_X^+)$, while the forgetful functor for $s$ sends $X$ to $(X,\mathcal{O}_X)$.
	That said, it should always be clear which forgetful functor is meant, as all the diagrams will only ever contain one of $p$ and $s$.
\end{remark}
How do the structure and support morphisms compare? One very basic observation takes place at the level of sets:
\begin{lemma}\label{psincl}
	Let $X$ be a locally Noetherian scheme, let $U$ be an open subscheme of $X$.
	Then we have an inclusion
	\[
	p^{-1}(U)\subseteq s^{-1}(U)
	\]
	in $X^{\ad}$.
\end{lemma}
\begin{proof}
	Let $\{V_i\}_i$ be an open affine cover of $X$.
	The restriction of $s_X$ to $V_i^{\ad}$ is $s_{V_i}$ since it is defined via gluing, and the restriction of $p_X$ to $V_i^{\ad}$ is $p_{V_i}$ by the definition of the adification.
	The $\{V_i^{\ad}\}_i$ form an open cover, so it is enough to show the claim for $X$ affine.
	
	In this case let $\nu\in\Spa{A}{A}$ be a valuation that is sent to $U$ under $p$, so $p(\nu)=\{\nu<1\}$ is in $U$.
	As $\{\nu=0\}$ is a generization of $\{\nu<1\}$ and $U$ is open, $s(\nu)=\{\nu=0\}$ is also contained in $U$.
\end{proof}
We have already explicitly seen what the pullback of $p$ to some open is in Corollary~\ref{schemestructurepb}, namely the adification of the open; however the pullback of $s$ is something different. Since this pullback will be very important in the proof of the theorem let us calculate a pullback of $s$ in the easiest possible case:
\begin{example}\label{distopen}
	Let $X=\Spec{A}$ be an Noetherian affine scheme and let $U\subseteq X$ be a distinguished open, say $D(f)$.
	We want to calculate $X\!\leftindex^{\ad\,}_s{\times}_XU$.
	
	The space $X^{\ad}=\Spa{A}{A}$ consists of all valuations of $A$ which are bounded by $1$.
	The map $s\colon X^{\ad}\to X$ sends a valuation to its support, so the preimage of $U$ consists of those valuations which do not map $f$ to $0$.
	Thus $X\!\leftindex^{\ad\,}_s{\times}_XU$ is the rational subspace $R(\frac{f}{f})$, which is given by
	\[
	\Spa{A\langle\frac{f}{f}\rangle}{A\langle\frac{f}{f}\rangle^+}=\Spa{A[f^{-1}]}{A^{\cl}}\mathrlap{,}
	\]
	where $(-)^{\cl}$ denotes integral closure.
	On the other hand the pullback of $U$ along $p$ is
	\[
	U^{\ad}=\Spa{A[f^{-1}]}{A[f^{-1}]}
	\]
	by Corollary~\ref{schemestructurepb}.
	Thus $X\!\leftindex^{\ad\,}_p{\times}_XU=U^{\ad}$ consists of those valuations in $X\!\leftindex^{\ad\,}_s{\times}_XU$ that send $f$ to $1$.
\end{example}
\begin{example}\label{distopenII}
	For a slightly more complicated example one can consider the preimage $W$ of $U$ under $s_X\circ i^{\ad}\colon\X^{\ad}\to X$, where $\X$ is the $f$-adic completion of~$X$.
	Similarly to the last example we will show that $W$ is given by
	\[
	R_{\X^{\ad}}(\frac{f}{f})=\Spa{A_f[f^{-1}]}{A^{\cl}}
	\]
	for some Huber ring $A_f[f^{-1}]$ which has $A_f$ as a ring of definition.
	
	We know that
	\[
	W\cong \U\times_{X^{\ad}}\X^{\ad}=\Spa{A[f^{-1}]}{A^{\cl}}\times_{\Spa{A}{A}}\Spa{A_f}{A}\mathrlap{.}
	\]
	For an adic space $T$ there exists an isomorphism
	\[
	\Ad(T,\Spa{A}{A^+})\cong\Huber((A,A^+),(\mathcal{O}_T(T),\mathcal{O}_T^+(T)))
	\]
	by \cite[Theorem 3.2.9.ii]{Hu}. Thus it is enough to show that $(A_f[f^{-1}],A^{\cl})$ is a pushout of
	\[
	\begin{tikzcd}
		(A,A)\arrow[r]\arrow[d]&(A_f,A)\\
		(A[f^{-1}],A^{\cl})
	\end{tikzcd}
	\]
	in the category of Huber pairs.
	Let us first define $(A_f[f^{-1}],A^{\cl})$ as a Huber pair.
	
	As a ring $A_f[f^{-1}]$ is given by $A[f^{-1}]$.
	Define a subset $U$ of $A_f[f^{-1}]$ which contains $0$ to be open iff it contains some $f^nA$.
	This yields a group topology by translating the open sets.
	We need to check that multiplication is continuous to obtain a topological ring structure on $A_f[f^{-1}]$.
	
	So let $x,y\in A_f[f^{-1}]$ and let $U$ be an open neighbourhood of $xy$.
	Then $U$ contains some $xy+f^nA$.
	Write $x=f^{-i}x^\prime$ and $y=f^{-j}y^\prime$ for $x^\prime, y^\prime\in A$.
	Then calculate
	\[
	(x+f^{j+n}A)(y+f^{i+n}A)=xy+(x^\prime+y^\prime)f^nA+f^{j+i+2n}A\subseteq xy+f^nA,
	\]
	showing continuity of the multiplication.
	The subring $A_f$ is clearly open, so $A_f[f^{-1}]$ is a Huber ring.
	Since $A$ is an open subring in $A_f[f^{-1}]$, so is $A^{\cl}$.
	But $A^{\cl}$ is additionally integrally closed, so $(A_f[f^{-1}],A^{\cl})$ is a Huber pair.
	
	The inclusion $A_f\to A_f[f^{-1}]$ is continuous, so we obtain a commutative square of maps of Huber pairs
	\[
	\begin{tikzcd}
		(A,A)\arrow[r]\arrow[d]&(A_f,A)\arrow[d]\\
		(A[f^{-1}],A^{\cl})\arrow[r]&(A_f[f^{-1}],A^{\cl})\mathrlap{.}
	\end{tikzcd}
	\]
	Let $(B,B^+)$ be a Huber pair, and let 
	\[
	\varphi\colon(A[f^{-1}],A^{\cl})\to(B,B^+)\text{ and } \psi\colon(A_f,A)\to(B,B^+)
	\]
	be maps of Huber pairs that agree on $(A,A)$.
	Since $(A,A)\to(A_f,A)$ is the identity on the underlying rings giving such a pair is the same as giving a morphism of Huber pairs
	\[
	\varphi\colon(A[f^{-1}],A^{\cl})\to(B,B^+)
	\]
	such that the restriction to $A$ is continuous for the $f$-adic topology.
	This means that for every neighbourhood of zero $U\subseteq B$ the set $A\,\cap\,\varphi^{-1}(U)$ is open in $A_f$, i.e.\@ it contains some $f^nA$.
	Due to $f^nA$ being contained in $A$ this is equivalent to every $\varphi^{-1}(U)$ containing $f^nA$, which means that $\varphi$ is continuous as a morphism $A_f[f^{-1}]\to B$.
	
	Thus we obtain a morphism of Huber pairs $(A_f[f^{-1}],A^{\cl})\to(B,B^+)$ which restricts to $\varphi$ and $\psi$.
	The map $(A[f^{-1}],A^{\cl})\to(A_f[f^{-1}],A^{\cl})$ is the identity on underlying rings, so this factorization is unique, showing the universal property.
\end{example}
Before we get to discussing coherent sheaves we need to note two more lemmata about adifications, fibre products and completions which will be important during the formulation of the main theorem later.
\begin{lemma}\label{structurecomm}
	Let $f\colon \widetilde{X}\to X$ be a morphism of locally Noetherian schemes and let $Z$ be a closed subscheme of $X$.
	Denote the completion of $X$ at $Z$ by $\X$ and the completion of $\widetilde{X}$ at $Z\times_X\widetilde{X}$ by $\widetilde{\X}$.
	Then the diagram
	\[
	\begin{tikzcd}
		\widetilde{\X}^{\ad}\arrow[r]\arrow[d]&\X^{\ad}\arrow[d]\\
		\widetilde{X}^{\ad}\arrow[r,"f^{\ad}"]\arrow[d,"s",swap]&X^{\ad}\arrow[d,"s"]\\
		\widetilde{X}\arrow[r,"f"]&X
	\end{tikzcd}
	\]
	commutes.
	If $f$ is locally of finite type it induces a map
	\[
	\zeta\colon\widetilde{\X}^{\ad}\to \X^{\ad}\times_X\widetilde{X}\mathrlap{,}
	\]
	which is an isomorphism if $f$ is proper.
\end{lemma}
\begin{proof}
	We just need to check the lower square for commutativity, as the upper square is the adification of a commutative square of formal schemes.
	Let $V$ be an open affine subscheme of $X$.
	The preimage of $V$ under $s$ contains $V^{\ad}$ by Lemma~\ref{psincl}.
	We can cover $X$ with affines $\{V_i\}_i$ and the preimages $f^{-1}(V_i)$ with affines $\{V^{\prime}_{ij}\}_{i,j}$.
	Restricting the maps to these covers yields a square
	\[
	\begin{tikzcd}
		V^{\prime\ad}_{ij}\arrow[r,"f^{\ad}"]\arrow[d,"s",swap]&V_i^{\ad}\arrow[d,"s"]\\
		V^{\prime}_{ij}\arrow[r,"f"]&V_i
	\end{tikzcd}
	\]
	of affine schemes and affinoid adic spaces. Let $V_i=\Spec{A}$ and $V^{\prime}_{ij}=\Spec B$.
	Both compositions
	\[
	(V^{\prime\ad}_{ij},\mathcal{O}_{V^{\prime\ad}_{ij}})\to (V_i,\mathcal{O}_{V_i})
	\]
	are determined by their map on global sections since $V_i$ is an affine scheme.
	The induced diagram on global sections is
	\[
	\begin{tikzcd}
		A\arrow[r,"f^\#"]\arrow[d,equal]&B\arrow[d,equal]\\
		A\arrow[r,"f^\#"]\arrow[r]&B\mathrlap{,}
	\end{tikzcd}
	\]
	since $s$ is induced by the identity by definition. Since this square commutes so does the square above.
	
	If $f$ is locally of finite type the space $\X^{\ad}\times_X\widetilde{X}$ exists by Lemma~\ref{fibreschemeadic}, so we obtain the map $\zeta$ by the universal property.
	For the addendum assume that $f$ is proper, in which case \cite[Example 3.10.6.v]{Hu} shows the claim.
\end{proof}
The following is the classical fact that the completion at a fibre product is the fibre product with the completion.
\begin{lemma}\label{formalpullback}
	Let $f\colon \widetilde{X}\to X$ be a morphism of locally Noetherian schemes and let $Z$ be a closed subscheme in $X$.
	Denote the completion of $X$ at $Z$ by $\X$ and denote the completion of $\widetilde{X}$ at $Z\times_X\widetilde{X}$ by $\widetilde{\X}$.
	Then the canonical morphism $\widetilde{\X}\to \X$ fits into a cartesian diagram of formal schemes
	\[
	\begin{tikzcd}
		\widetilde{\X}\arrow[r]\arrow[d]&\X\arrow[d]\\
		\widetilde{X}\arrow[r]&X\mathrlap{.}
	\end{tikzcd}
	\]
\end{lemma}
\begin{proof}
	This is \cite[I.10.9.9]{EGA}.
\end{proof}
\subsection{Coherent sheaves}
To define what coherent sheaves are one can try to emulate the situation for schemes, where there is a local condition for coherence. Thus one needs to first define which kind of sheaf the restriction of a coherent sheaf to some affinoid should be.
When following Hubers exposition in \cite[Chapter 3.6]{Hu} one needs the adic space to locally have Noetherian rings of definition to be able to define a good notion of coherent sheaves.
Thus we will only define coherence for locally Noetherian adic spaces, in which case it will turn out that the following definition of coherence agrees with the usual definition of a coherent $\mathcal{O}_{(X,\mathcal{O}_X)}$-module.
Since all adic spaces that we encounter are locally Noetherian by Lemma~\ref{Noetherian} it will be enough to consider this case.

Until the end of the section let $X$ be a locally Noetherian adic space.
\begin{definition}\label{cohdef}
	Let $X=\Spa{A}{A^+}$ be affinoid. Let $M$ be a finitely generated $\mathcal{O}_X(X)$-module.
	For a rational open $V\subseteq X$ define
	\[
	(M\otimes\mathcal{O}_X)(V):=M\otimes_{\mathcal{O}_X(X)} \mathcal{O}_X(V)\mathrlap{.}
	\]
	For a general open subset $W\subseteq X$ we can now define
	\[
	(M\otimes\mathcal{O}_X)(W):=\lim_{V\subseteq W}(M\otimes\mathcal{O}_X)(V)\mathrlap{,}
	\]
	where $V$ runs through all rational subsets of $X$ which are contained in $W$.
	This yields an $\mathcal{O}_X$-module presheaf $M\otimes\mathcal{O}_X$.
\end{definition}
\begin{proposition}
	Let $X=\Spa{A}{A^+}$ be affinoid and let $M$ be a finitely generated $\mathcal{O}_X(X)$-module.
	Then $M\otimes\mathcal{O}_X$ is a sheaf.
\end{proposition}
Before we prove this proposition, let us remark on the difference of this definition to the definition in \cite{Hu} and the literature more generally.
\begin{remark}
	In \cite[Chapter 3.6]{Hu} the sheaf $M\otimes\mathcal{O}_X$ is defined the same way as above, but one additionally equips $(M\otimes\mathcal{O}_X)(V)$ with the quotient topology coming from writing it as a quotient of some $\mathcal{O}_X(V)^n$. Then the limit is equipped with the limit topology, yielding a presheaf of complete topological modules.
	This is the natural thing to take when trying to talk about quasicoherent sheaves on adic spaces, but as we will see later that it does not make a difference whether one considers coherent sheaves of topological modules or abstract modules.
\end{remark}
\begin{proof}
	Huber shows that $M\otimes\mathcal{O}_X$ equipped with this natural topology is a sheaf of topological modules in \cite[Theorem 3.6.2]{Hu}. As the forgetful functor from topological modules to abstract modules is continuous the proposition follows.
\end{proof}
\begin{definition}
	An $\mathcal{O}_X$-module $\mathcal{M}$ is called \emph{coherent} if for every affinoid open $V\subseteq X$ there exists a finitely generated $\mathcal{O}_X(V)$-module $M$ such that
	\[
	M\otimes\mathcal{O}_V\cong\mathcal{M}|_V\mathrlap{.}
	\]
	The full subcategory of the category of $\mathcal{O}_X$-modules consisting of coherent $\mathcal{O}_X$-modules is denoted by $\Coh[X]$.
\end{definition}
\begin{remark}\item
	\begin{itemize}
		\item 
		When using Huber's definition of $M\otimes\mathcal{O}_X$ the definition of coherent modules is completely analogous as a full subcategory of the category of topological $\mathcal{O}_X$-modules. We will call this category $\Coht[X]$, the category of coherent topological sheaves.
		\item 
		In the definition one can equivalently require that there exists an open affinoid cover over which $\mathcal{M}$ is given by a module $M$.
	\end{itemize}
\end{remark}
We note that coherent (topological) sheaves satisfy descent for open covers:
\begin{lemma}\label{zardesc}
	Let $\{V_i\}_i$ be an open cover of $X$, then coherent (topological) sheaves on $X$ satisfy descent for $\{V_i\}_i$.
\end{lemma}
\begin{proof}
	By Lemma~\ref{pullback} we have $V_i\times_X V_j=V_i\,\cap\, V_j$. Over opens and their intersections we are able to glue morphisms uniquely as with modules on any ringed space.
	Moreover all descent data for (topological) $\mathcal{O}_X$-modules are effective; since being coherent can be checked on an open cover they are also effective for $\Coh[X]$ and $\Coht[X]$.
\end{proof}

Let us now check that it does indeed not make a difference whether our coherent sheaves are equipped with a topology or not.
\begin{lemma}
	The forgetful functor $\Coht[X]\to\Coh[X]$ is an equivalence.
\end{lemma}
\begin{proof}
	First let us note that this functor exists by continuity of the forgetful functor from topological modules to modules.
	It is enough to check the claim for $X$ affinoid, since both sides satisfy descent for open covers. Essential surjectivity and faithfulness are clear.
	
	Let $M$ and $N$ be finitely generated $\mathcal{O}_X$-modules.
	For fullness it is enough to show that any morphism
	\[
	\begin{tikzcd}
		M\otimes_{\mathcal{O}_X(X)} \mathcal{O}_X(V)\arrow[r,"f"] & N\otimes_{\mathcal{O}_X(X)} \mathcal{O}_X(V)
	\end{tikzcd}
	\]
	with $V\subseteq X$ rational is continuous. Taking generators of $M$ and $N$ yields a commutative diagram
	\[
	\begin{tikzcd}
		&\mathcal{O}_X(V)^n\arrow[dl,two heads]\arrow[dr,two heads]\\
		M\otimes_{\mathcal{O}_X(X)} \mathcal{O}_X(V) \arrow[rr,"f"]&& N\otimes_{\mathcal{O}_X(X)} \mathcal{O}_X(V)\mathrlap{,}
	\end{tikzcd}
	\]
	where the diagonal morphisms are topological quotient maps, so $f$ is continuous.
\end{proof}
Thus we are able to work with $\Coh[X]$ for the remainder of the text, even when using results from \cite{Hu}. As it turns out this category is just the usual category of coherent sheaves.
\begin{lemma}\label{cohrspace}
	The category $\Coh[X]$ is the category of coherent sheaves of $(X,\mathcal{O}_X)$ considered as a ringed space.
\end{lemma}
\begin{proof}
	Both of these are full subcategories of the category of $\mathcal{O}_X$-modules, so it is enough to check that their objects agree. Both satisfy descent for open covers, so we can assume that $X$ is affinoid.
	Since $X$ is a locally Noetherian adic space $\mathcal{O}_X(V)$ is Noetherian for every rational $V\subseteq X$ by \cite[Proposition 3.6.5.ii]{Hu}. Thus a module on $(X,\mathcal{O}_X)$ is coherent iff it is quasicoherent and of finite type.
	Any quasicoherent and finite type module on a ringed space is the sheafification of a presheaf $\mathcal{F}(U)=M\otimes_{\mathcal{O}_X(X)}\mathcal{O}_X(U)$ for some finitely generated $\mathcal{O}_X(X)$-module $M$.
	Thus it is enough to check that this sheafification is $M\otimes\mathcal{O}_X$.
	It is a short verification that the map $\mathcal{F}\to M\otimes\mathcal{O}_X$ which sends a section on $U$ to all its restrictions fulfills the universal property of the sheafification.
\end{proof}
\begin{proposition}\label{globalsec}
	Let $X=\Spa{A}{A^+}$ be affinoid.
	Then the global sections functor
	\[
	\Gamma\colon\Coh[X]\to \fMod[\mathcal{O}_X(X)]
	\]
	is an equivalence of categories. There exists an exact inverse functor $-\otimes\mathcal{O}_X$ which sends a module $M$ to $M\otimes\mathcal{O}_X$.
\end{proposition}
\begin{proof}
	Combine \cite[Proposition 3.6.19.iii]{Hu} and \cite[Theorem 3.6.20]{Hu}.
\end{proof}
\subsection{Pullbacks of coherent sheaves}
Considering Lemma~\ref{cohrspace} we are able to pull back coherent sheaves via any morphism of adic spaces as well as $s_X$ by viewing them as morphisms of ringed spaces.
Coherence is preserved by these pullbacks since all spaces are locally Noetherian (as ringed spaces), so being coherent is the same as being quasicoherent and of finite type.
\begin{remark}\item
	\begin{itemize}
		\item 
		If $\iota\colon U\to X$ is an open immersion of locally Noetherian adic spaces then $\iota^*$ identifies with the restriction functor $\Coh[X]\to\Coh[U]$.
		\item 
		The functor $-\otimes\mathcal{O}_X$ is the pullback on coherent sheaves along
		\[
		\pi\colon(X,\mathcal{O}_X)\to(\pt,\mathcal{O}_X(X))\mathrlap{,}
		\]
		where the map on structure sheaves is the identity on $\mathcal{O}_X(X)$.
		\item 
		Later on we will also need to pull back coherent sheaves along the structure map
		\[
		(\X^{\ad},\mathcal{O}^+_{{\X}^{\ad}})\xrightarrow{p}(\X,\mathcal{O}_{X})
		\]
		of an adification, by which we will mean that we pull back along the composition
		\[
		(\X^{\ad},\mathcal{O}_{\X^{\ad}})\to(\X^{\ad},\mathcal{O}^+_{{\X}^{\ad}})\xrightarrow{p}(\X,\mathcal{O}_{X})\mathrlap{.}
		\]
	\end{itemize}
\end{remark}
For morphisms between affinoid spaces there is an explicit description of the pullback, which should be familiar from scheme theory. On global sections the pullback is given by tensoring with the global sections:
\begin{lemma}\label{pullbackaffinoid}
	Let each of $X$ and $Y$ be either an affine formal scheme or an affinoid adic space, and let $f\colon (X,\mathcal{O}_X)\to(Y,\mathcal{O}_Y)$ be a morphism of ringed spaces.
	Then the following diagram commutes:
	\[
	\begin{tikzcd}[column sep=60]
		\Coh[Y]\arrow[r,"f^*"]\arrow[d,"\Gamma",swap]\arrow[d,"\cong"]&\Coh[X]\arrow[d,"\Gamma"]\arrow[d,"\cong",swap]\\
		\fMod[\mathcal{O}_Y(Y)]\arrow[r,"-\otimes_{\mathcal{O}_Y(Y)}\mathcal{O}_X(X)"]&\fMod[\mathcal{O}_X(X)]
	\end{tikzcd}
	\]
\end{lemma}
\begin{proof}
	Let $f^\#$ be the morphism $f$ induces on global sections.
	Consider the following commutative diagram of ringed spaces:
	\[
	\begin{tikzcd}
		(X,\mathcal{O}_X)\arrow[r,"f"]\arrow[d,"\pi",swap]&(Y,\mathcal{O}_Y)\arrow[d,"\pi"]\\
		(\pt,\mathcal{O}_X(X))\arrow[r,"f^\#"]&(\pt,\mathcal{O}_Y(Y))
	\end{tikzcd}
	\]
	The commutative square induced by pulling back is
	\[
	\begin{tikzcd}[column sep=60]
		\fMod[\mathcal{O}_Y(Y)]\arrow[d,"-\otimes\mathcal{O}_Y",swap]\arrow[r,"-\otimes_{\mathcal{O}_Y(Y)}\mathcal{O}_X(X)"]&\fMod[\mathcal{O}_X(X)]\arrow[d,"-\otimes\mathcal{O}_X"]\\
		\Coh[Y]\arrow[r,"f^*"]&\Coh[X]\mathrlap{,}
	\end{tikzcd}
	\]
	where for an affine formal scheme $X$ the functor $-\otimes\mathcal{O}_X$ sends a module to the associated coherent sheaf on $X$ (this is often called $\widetilde{(-)}$ instead; the above notation treats both cases simultaneously).
	Thus the square of the lemma also commutes.
\end{proof}
\begin{proposition}\label{isoU}
	Let $X$ be a locally Noetherian scheme, $U$ an open subscheme. Let $\U$ be the preimage of $U$ under $s\colon X^{\ad}\to X$.
	Then the inclusion $U^{\ad}\to \U$ induces an equivalence on coherent sheaves.
\end{proposition}
\begin{proof}
	We can cover $X^{\ad}$ by $V^{\ad}$ for $V\subseteq X$ affine, so we can cover $\U$ by
	\[
	\U\times_{X^{\ad}}V^{\ad}=(U\times_X X^{\ad})\times_{X^{\ad}}V^{\ad}\cong U\times_X V^{\ad}\cong 	(U\times_XV)\times_VV^{\ad}\mathrlap{.}
	\]
	The intersection of this with $U^{\ad}$ is
	\[
	U^{\ad}\times_{X^{\ad}}\U\times_{X^{\ad}}V^{\ad}=U^{\ad}\times_{X^{\ad}}V^{\ad}\cong (U\times_XV)^{\ad}
	\]
	by Lemma~\ref{openint}.
	Thus using descent along open covers from Lemma~\ref{zardesc} we can assume that $X=\Spec A$ is affine and that $U=D(f)$ is a distinguished open.
	In this case we have 
	\[
	U^{\ad}=\Spa{A[f^{-1}]}{A[f^{-1}]}
	\]
	and
	\[
	\U\cong\Spa{A[f^{-1}]}{A^{\cl}}
	\]
	by Corollary~\ref{schemestructurepb} and Example~\ref{distopen}.
	By the construction of $\mathcal{O}_X$ the map between them is induced by the inclusion
	\[
	(A[f^{-1}],A^{\cl})\to(A[f^{-1}],A[f^{-1}])
	\]
	of Huber pairs, which is the identity on global sections. 
	Using Lemma~\ref{pullbackaffinoid} and Proposition~\ref{globalsec} we conclude that the pullback functor
	\[
	\Coh[\U]\to\Coh[U^{\ad}]
	\]
	is an equivalence.
\end{proof}
In the proof of the theorem it will be important to relate coherent sheaves on a formal scheme to coherent sheaves on its adification.
\begin{lemma}\label{pequiv}
	Let $\X$ be a locally Noetherian formal scheme.
	The pullback functor
	\[
	p^*\colon \Coh[\X]\to\Coh[\X^{\ad}]
	\]
	is an equivalence.
\end{lemma}
\begin{proof}
	We can assume that $\X=\Spf{A}$ is affine, since the preimage of an open formal subscheme $U\subseteq \X$ under $p$ is $U^{\ad}$, adification commutes with intersection and coherent sheaves have descent for open covers.
	In this case both the global sections of $\mathcal{O}_{\X^{\ad}}$ and $\mathcal{O}^+_{\X^{\ad}}$ are equal to $\comp{A}$, and the morphism of ringed spaces
	\[
	(\X^{\ad},\mathcal{O}_{\X}^{\ad})\to(\X^{\ad},\mathcal{O}^+_{\X^{\ad}})
	\]
	induces the identity on global sections.
	Moreover by Remark~\ref{pglobalsec} the structure map $p$ is the preimage of the identity on global sections under the isomorphism
	\[
	\TLRS((\X^{\ad},\mathcal{O}^+_{{\X}^{\ad}}),(\Spf{A},\mathcal{O}_{\Spf{A}}))\to 	\TopRing(\mathcal{O}_{\Spf{A}}(\Spf{A}),\mathcal{O}^+_{\X^{\ad}}(\X^{\ad}))\mathrlap{,}
	\]
	which sends a morphism to its action on global sections.
	Thus the following diagram commutes
	\[
	\begin{tikzcd}
		\Coh[\X]\arrow[r]\arrow[d,"\Gamma",swap]\arrow[d,"\cong"]&\Coh[(\X^{\ad},\mathcal{O}^+_{{\X}^{\ad}})]\arrow[d,"\Gamma"]\arrow[r]&\Coh[\X^{\ad}]\arrow[d,"\Gamma"]\arrow[d,"\cong",swap]\\
		\fMod[\comp{A}]\arrow[r,equal]&\fMod[\comp{A}]\arrow[r,equal]&\fMod[\comp{A}]
	\end{tikzcd}
	\]
	by Lemma~\ref{pullbackaffinoid}, so $p^*$ is an equivalence.
\end{proof}
The next lemma enables us through Corollary~\ref{pushforward_iso} to build the commutative diagram setting up the theorem in Section~\ref{setup}, in which the pullback for a scheme $X$ should be taken along $s$, and the pullback for its completion $\X$ has to be taken along $p$, as there is no morphism $s$ for formal schemes.
\begin{lemma}\label{pequals}
	Let $X$ be a locally Noetherian scheme.
	The pullback functors $s^*$ and $p^*$ are equal.
\end{lemma}
\begin{proof}
	Let $X$ be a locally Noetherian scheme, and let $U$ be an open subscheme.
	The restrictions of both $s_X$ and $p_X$ to $U^{\ad}\to U$ are $s_U$ and $p_U$, as $s$ and $p$ are defined by gluing locally.
	Since coherent sheaves have descent for open covers we can thus assume $X$ to be affine.
	
	In this case both $s$ and $p$ induce the identity on global sections by definition. Moreover
	\[
	(X^{\ad},\mathcal{O}_{X^{\ad}})\to(X^{\ad},\mathcal{O}^+_{X^{\ad}})
	\]
	induces the identity on global sections, so by Lemma~\ref{pullbackaffinoid} the pullbacks agree.
\end{proof}
\begin{corollary}\label{sequiv}
	Let $X$ be a locally Noetherian scheme, then $s^*$ is an equivalence.
\end{corollary}
\begin{proof}
	This follows from Lemma~\ref{pequals} and Lemma~\ref{pequiv}.
\end{proof}
\begin{corollary}\label{pushforward_iso}
	Let $X$ be a locally Noetherian scheme, let $Z$ be a closed subscheme of $X$. Denote the completion of $X$ at $Z$ by $\X$.
	Then the following square commutes:
	\[
	\begin{tikzcd}
		\Coh[X]\arrow[r,"s^*"]\arrow[d,"i^*",swap]&\Coh[X^{\ad}]\arrow[d,"i^{\ad,*}"]\\
		\Coh[\X]\arrow[r,"p^*"]&\Coh[\X^{\ad}]
	\end{tikzcd}
	\]
\end{corollary}
\begin{proof}
	The diagram
	\[
	\begin{tikzcd}
		(\X^{\ad},\mathcal{O}_{\X^{\ad}})\arrow[r]\arrow[d,"i^{\ad}",swap]&(\X^{\ad},\mathcal{O}^+_{\X^{\ad}})\arrow[r,"p_{\X}"]\arrow[d,"i^{\ad}"]&(\X,\mathcal{O}_{\X})\arrow[d,"i"]\\
		(X^{\ad},\mathcal{O}_{X^{\ad}})\arrow[r]&(X^{\ad},\mathcal{O}^+_{X^{\ad}})\arrow[r,"p_X"]&(X,\mathcal{O}_X)
	\end{tikzcd}
	\]
	clearly commutes, so the corresponding square of pullback functors commutes as well. Then the claim follows by Lemma~\ref{pequals}.
\end{proof}
\newpage
\section{Setting up the theorem}\label{setup}
In this section we want to explain how the setup of the theorem works,  and check that it plays nicely with basechange.
We fix some notation for this section:
\begin{itemize}
	\item 
	$X$ a locally Noetherian scheme with a closed subscheme $Z$
	\item 
	$i\colon\X\to X$ the completion of $X$ at $Z$
	\item 
	$j\colon U\to X$ the complement of $Z$ in $X$
	\item 
	$\U$ the preimage of $U$ under $s\colon X^{\ad}\to X$
	\item 
	$W$ the pullback of $\U\to X^{\ad}$ along $i^{\ad}$
\end{itemize}
Here the pullback $W$ exists and is given by the preimage of $\U$ in $\X^{\ad}$ by Lemma~\ref{pullback}.
This means equivalently that $W$ is the preimage of $U$ under $s\circ i^{\ad}$, since $\U$ is the preimage of $U$ under $s$.
We can summarize the notation above in a commutative diagram
\[
\begin{tikzcd}
	W\arrow[r]\arrow[d]&\X^{\ad}\arrow[d,"i^{\ad}"]\\
	\U\arrow[r]\arrow[d,"s",swap]&X^{\ad}\arrow[d,"s"]\\
	U\arrow[r,"j"]&X\mathrlap{.}
\end{tikzcd}
\]
Pulling back along these maps and using Corollary~\ref{pushforward_iso} yields a commutative diagram
\[
\begin{tikzcd}
	&\Coh[X]\arrow[r,"j^*"]\arrow[d,"s^*"]\arrow[ddl,"i^*",bend right=20,swap]&\Coh[U]\arrow[d,"s^*"]\\
	&\Coh[X^{\ad}]\arrow[r]\arrow[d,"i^{\ad,*}"]&\Coh[\U]\arrow[d]\\
	\Coh[\X]\arrow[r,"p^*"]&\Coh[\X^{\ad}]\arrow[r]&\Coh[W]\mathrlap{.}
\end{tikzcd}
\]
When referencing the ``setup of the theorem for $(X,Z)$'' we mean this diagram, while the outer square will be called the ``square of the theorem for $(X,Z)$''. 
Now we can precisely formulate the descent theorem that we want to prove:
\begin{theorem}[Main theorem]\label{mainthmprop}
	Let $X$ be a locally Noetherian scheme and let $Z$ be a closed subscheme of $X$.
	In the setup above the outer square
	\[
	\begin{tikzcd}
		\Coh[X]\arrow[r,"j^*"]\arrow[d,"i^*",swap]&\Coh[U]\arrow[d]\\
		\Coh[\X]\arrow[r]&\Coh[W]
	\end{tikzcd}
	\]
	is cartesian.
\end{theorem}
Let us consider the way in which the setup is functorial; for the rest of the section fix the following additional notation.
Let $f\colon\widetilde{X}\to X$ be a morphism of locally Noetherian schemes, denote pullbacks of the objects in the setup above along $f$ by a $\widetilde{(-)}$.
Consider the following diagram
\[
\begin{tikzcd}
	\widetilde{\X}^{\ad}\arrow[r]\arrow[d]&\X^{\ad}\arrow[d]\\
	\widetilde{X}^{\ad}\arrow[r,"f^{\ad}"]\arrow[d,"s",swap]&X^{\ad}\arrow[d,"s"]\\
	\bX\arrow[r,"f"]&X\mathrlap{,}
\end{tikzcd}
\]
which is commutative by Lemma~\ref{structurecomm}.

The pullback of this diagram to $U$ is the commutative diagram
\[
\begin{tikzcd}
	\bW\arrow[r]\arrow[d]&W\arrow[d]\\
	\widetilde{\U}\arrow[r]\arrow[d,"s",swap]&\U\arrow[d,"s"]\\
	\widetilde{U}\arrow[r]&U\mathrlap{,}
\end{tikzcd}
\]
since $f^{-1}(U)=\widetilde{U}$.
Applying $\Coh[-]$ and Corollary~\ref{pushforward_iso} and pasting the resulting diagrams one obtains a large commutative diagram
\[
\begin{tikzcd}
	&\Coh[X]\arrow[rr,"j^*"]\arrow[dd,"s^*"]\arrow[ddddl,"i^*",swap]\arrow[dr]&&\Coh[U]\arrow[dd,"s^*",near start]\arrow[dr]\\
	&&\Coh[\bX]\arrow[rr,crossing over]&&\Coh[\widetilde{U}]\arrow[dd,"s^*"]\\
	&\Coh[X^{\ad}]\arrow[rr]\arrow[dd,"i^{\ad,*}"]\arrow[dr]&&\Coh[\U]\arrow[dd]\arrow[dr]\\
	&&\Coh[\widetilde{X}^{\ad}]\arrow[rr,crossing over]\arrow[from=uu,crossing over,"s^*",near start]&&\Coh[\widetilde{\U}]\arrow[dd]\\
	\Coh[\X]\arrow[dr]\arrow[r,"p^*"]&\Coh[\X^{\ad}]\arrow[rr]\arrow[dr]&&\Coh[W]\arrow[dr]\\
	&\Coh[\widetilde{\X}]\arrow[r]\arrow[from=uuuur,end anchor={[xshift=0.8em]},crossing over]&\Coh[\widetilde{\X}^{\ad}]\arrow[from=uu,crossing over]\arrow[rr]&&\Coh[\bW]\mathrlap{.}
\end{tikzcd}
\]
The back face of the large diagram is the setup of the theorem for $(X,Z)$ while the front face is the setup of the theorem for $(\bX,Z^\prime)$.
The diagonal functors are given by pullback along the appropriate basechange of $f\colon \bX\to X$.

We can summarize the discussion overhead in a lemma:
\begin{lemma}\label{functsetup}
	With the notation above there exists a commutative cube
	\[
	\begin{tikzcd}
		\Coh[X]\arrow[rr]\arrow[dr]\arrow[dd]&&\Coh[U]\arrow[dd]\arrow[dr]\\
		&\Coh[\bX]\arrow[rr,crossing over]&&\Coh[\widetilde{U}]\arrow[dd]\\
		\Coh[\X]\arrow[rr]\arrow[dr]&&\Coh[W]\arrow[dr]\\
		&\Coh[\widetilde{\X}]\arrow[rr]\arrow[from=uu,crossing over]&&\Coh[\bW]
	\end{tikzcd}
	\]
	in which the back face is the square of the theorem for $(X,Z)$, the front face is the square of the theorem for $(\bX,Z^\prime)$ and the diagonal functors are pullbacks along an appropriate basechange of $f\colon \bX\to X$.\qed
\end{lemma}
\begin{remark}\label{pbsquare}
	\item 
	\begin{itemize}
		\item 
		The reason why gluing along $W$ yields a different result from the naive attempt at pulling back on schemes is the map we are pulling back by.
		The preimage of $U$ under $p\circ i^{\ad}$ is the adification of $\X\times_XU$ by Corollary~\ref{structurepreimage}, which is empty.
		When pulling back by $s\circ i^{\ad}$ instead the preimage $W$ is a slightly bigger subspace (as f.e.\@ evidenced by Lemma~\ref{psincl}), which turns out to be large enough to obtain a descent statement.
		\item 
		Define
		\[
		\C_{X}\coloneqq\Coh[\X]\times_{\Coh[W]}\Coh[U]
		\]
		and denote the pullback functor of the theorem by
		\[
		\Phi_X\colon \Coh[X]\to\C_X\mathrlap{.}
		\]
		Then the commutative cube above becomes a commutative square
		\[
		\begin{tikzcd}
			\Coh[\bX]\arrow[d]\arrow[r,"\Phi_{\bX}"]&\C_{\bX}\arrow[d]\\
			\Coh[X]\arrow[r,"\Phi_X"]&\C_X\mathrlap{.}
		\end{tikzcd}
		\]
		In particular we can apply this to open subschemes $V\subseteq \widetilde{V}$ in $X$.
		Applying the commutative cube in Lemma~\ref{functsetup} to $V$ as a subscheme of $\widetilde{V}$ yields a commutative square
		\[
		\begin{tikzcd}
			\Coh[\widetilde{V}]\arrow[d,"\iota^*",swap]\arrow[r,"\Phi_{\widetilde{V}}"]&\C_{\widetilde{V}}\arrow[d]\\
			\Coh[V]\arrow[r,"\Phi_V"]&\C_V\mathrlap{,}
		\end{tikzcd}
		\]
		where the functor $\C_{\widetilde{V}}\to \C_V$ is given by the respective pullbacks along the restrictions.
	\end{itemize}
\end{remark}
Now that we have seen how the setup of the theorem is functorial the following lemma will allow us to reduce to the affine case during the proof of the theorem.
\begin{lemma}\label{affred}
	Assume that for every open subscheme $V\subseteq X$ there exists an open cover $\{V_i\}_i$ of $V$ such that the theorem holds for every $(V_i,Z\cap V_i)$.
	Then the theorem also holds for $(X,Z)$.
\end{lemma}
\begin{proof}
	First note that for an open cover $\{V_i\}_i$ of $X$ the collections $\{V_i\times_X\X\}_i$, $\{V_i\times_XU\}_i$ and $\{V_i^{\ad}\times_{X^{\ad}}W\}_i$ appearing when applying Lemma~\ref{functsetup} to $\coprod V_i\to X$ are open covers of $\X$, $U$ and $W$, respectively.
	
	By the assumption we can choose an open cover $\{V_i\}_i$ of $X$ such that the theorem holds for all $(V_i,Z\,\cap\, V_i)$.
	We can also choose an open cover $\{V^\prime_{ijk}\}_k$ of $V_{ij}\coloneqq V_i\,\cap\, V_j$ such that the theorem holds for all $(V^\prime_{ijk},Z\,\cap\, V^\prime_{ijk})$.
	Repeatedly pasting the commutative square from Remark~\ref{pbsquare} for the different inclusions between the subschemes we have just chosen yields a commutative diagram
	\[
	\begin{tikzcd}[column sep=50]
		\Coh[X]\arrow[r,"\Phi_X"]\arrow[d]&\C_X\arrow[d]\\
		\prod_i\Coh[V_i]\arrow[r,"\cong",swap]\arrow[r,"(\Phi_i)_{i}"]\arrow[d,shift left]\arrow[d,shift right]
		&\prod_i\C_{V_i}\arrow[d,shift left]\arrow[d,shift right]\\
		\prod_{i,j}\Coh[V_{ij}]\arrow[r,"(\Phi_{ij})_{i,j}"]\arrow[d]
		&\prod_{i,j}\C_{V_{ij}}\arrow[d]\\
		\prod_{i,j,k}\Coh[V^\prime_{ijk}]\arrow[r,"\cong",swap]\arrow[r,"(\Phi_{ijk})_{i,j,k}"]&\prod_{i,j,k}\C_{V^\prime_{ijk}}\mathrlap{,}
	\end{tikzcd}
	\]
	where the single vertical functors are restrictions while the parallel functors are given by the projection to the $i$-th factor and restriction, or the projection to the $j$-th factor and restriction, respectively.
	The vertical single functors are given by restricting sheaves onto open covers, so they are all faithful.
	Thus both $\Phi_X$ and $(\Phi_{ij})_{i,j}$ are faithful.
	Let $F$ and $G$ be two coherent sheaves on $X$. Since $\{V_i\}_i$ is an open cover of $X$ and the pullbacks to $\X$, $U$ and $W$ are open covers as well the two columns induce equalizer diagrams
	\[
	\begin{tikzcd}
		0\arrow[d]&0\arrow[d]\\
		\Coh[X](F,G)\arrow[d]\arrow[r,hook]&\C_X(\Phi_X(F),\Phi_X(G))\arrow[d]\\
		\prod_i\Coh[V_i](F|_{V_i},G|_{V_i})\arrow[d,shift left]\arrow[d,shift right]\arrow[r,"\cong"]
		&\prod_i\C_{V_i}(\Phi_X(F)|_{V_i},\Phi_X(G)|_{V_i})\arrow[d,shift left]\arrow[d,shift right]\\
		\prod_{i,j}\Coh[V_{ij}](F|_{V_{ij}},G|_{V_{ij}})\arrow[r,hook]&\prod_{i,j}\C_{V_{ij}}(\Phi_X(F)|_{V_{ij}},\Phi_X(G)|_{V_{ij}})\mathrlap{.}
	\end{tikzcd}
	\]
	Both $\Phi_X$ and $(\Phi_{ij})_{i,j}$ are faithful and $(\Phi_i)_i$ is an equivalence, so the respective morphisms above are injective or an isomorphism. Thus the map
	\[
	\Coh[X](F,G)\to\C_X(\Phi_X(F),\Phi_X(G))
	\]
	is an isomorphism by a short diagram chase, i.e.\@ $\Phi_X$ is full. Redoing the entire argument for $X=V_{ij}$ one sees that $(\Phi_{ij})_{i,j}$ is fully faithful as well.
	
	Now we need to show essential surjectivity of $\Phi_X$. By Lemma~\ref{zardesc} a coherent sheaf on $X$ is equivalently given by coherent sheaves $M_i$ on $X_i$ together with isomorphisms
	\[
	\varphi_{ij}\colon M_j|_{X_{ij}}\to M_i|_{X_{ij}}
	\]
	that fulfill the cocycle condition.
	Similarly an object $N$ of $\C_X$ is equivalently given by objects $N_i$ in $\C_{V_i}$ together with isomorphisms
	\[
	\psi_{ij}\colon N_j|_{V_{ij}}\to N_i|_{V_{ij}}
	\]
	fulfilling the cocycle condition.
	Thus we can alternatively check that the induced functor between the categories of descent data is essentially surjective.
	
	Since $(\Phi_i)_i$ is an equivalence and $(\Phi_{ij})_{i,j}$ is fully faithful we can lift the sheaves and the isomorphisms of any descent datum $(N_i,\psi_{ij})$ to a collection $(M_i,\varphi_{ij})$ up to isomorphism. If the $\varphi_{ij}$ fulfill the cocycle condition this collection is a descent datum on $\Coh[X]$ mapping to $(N_i,\psi_{ij})$ up to isomorphism, so let us check the cocycle condition.
	The functor
	\[
	\Phi_{V_{ijk}}\colon\Coh[V_{ijk}]\to\C_{V_{ijk}}
	\]
	on triple intersections is faithful by the first half of the proof.
	The images under $\Phi_{V_{ijk}}$ of the restrictions of the $\varphi_{ij}$ to $V_{ijk}$ are the restrictions of the $\psi_{ij}$ to $V_{ijk}$. These fulfill the cocycle condition, so by full faithfulness the $\varphi_{ij}$ do so as well.
	Thus the functor on descent data is essentially surjective and so is $\Phi_X$.
\end{proof}
Later we want to reduce to the case of $Z$ being locally cut out by a principal ideal in the proof by blowing up $X$ in $Z$. For this we need to show that $\U$ and $W$ are invariant under blowups $Z$.
\begin{lemma}\label{modification}
	Assume that $f\colon \bX\to X$ is proper and an isomorphism over $U$.
	Then the maps $\widetilde{\U}\to \U$ and $\bW\to W$ used in Lemma~\ref{functsetup} are isomorphisms of adic spaces.
\end{lemma}
\begin{proof}
	Recalling the setup the maps $\widetilde{\U}\to\U$ and $\bW\to W$ are given as the pullback along $U$ of the maps $f^{\ad}\colon \widetilde{X}^{\ad}\to X^{\ad}$ and $\widetilde{f}^{\ad}\colon \widetilde{\X}^{\ad}\to \X^{\ad}$, so they fit into a commutative cube
	\[
	\begin{tikzcd}
		\bW\arrow[rr]\arrow[dr]\arrow[dd]&&W\arrow[dd]\arrow[dr]\\
		&\widetilde{\X}^{\ad}\arrow[rr,crossing over]&&\X^{\ad}\arrow[dd]\\
		\widetilde{\U}\arrow[rr]\arrow[dr]&&\U\arrow[dr]\\
		&\widetilde{X}^{\ad}\arrow[rr]\arrow[from=uu, crossing over]&&X^{\ad}\mathrlap{.}
	\end{tikzcd}
	\]
	Since $f$ is proper the support map and $f^{\ad}$ induce an isomorphism
	\[
	\widetilde{X}^{\ad}\to X^{\ad}\times_X\bX
	\]
	by Lemma~\ref{structurecomm}.
	Using this map we obtain the following chain of isomorphisms
	\[
	\begin{aligned}
		\widetilde{\U}&=\widetilde{\X}^{\ad}\times_{\bX}\widetilde{U}\\
		&\cong(X^{\ad}\times_X\bX)\times_{\bX}\widetilde{U}\\
		&\cong X^{\ad}\times_X\widetilde{U}\\
		&\cong X^{\ad}\times_XU=\U
	\end{aligned}
	\]
	where the last step uses the assumption that $f$ is an isomorphism over $U$.
	This composition is given by $f^{\ad}$ on the first component and the restriction of $f$ on the second, so it is the map we are interested in.
	Similarly we obtain a map
	\[
	\widetilde{\X}^{\ad}\to \X^{\ad}\times_X\bX
	\]
	which, using properness of $f$, is also an isomorphism by Lemma~\ref{structurecomm}.
	Moreover we know 
	\[
	W=\X^{\ad}\times_{X^{\ad}}\U\cong \X^{\ad}\times_XU\
	\]
	and
	\[
	\bW\cong \widetilde{\X}^{\ad}\times_XU\mathrlap{,}
	\]
	so we can write $\bW\to W$ as a similar chain of isomorphisms.
\end{proof}
\begin{corollary}\label{blowupsetup}
	Assume that $f\colon \bX\to X$ is proper and an isomorphism over $U$.
	Combining Lemma~\ref{functsetup} and Lemma~\ref{modification} yields a commutative diagram
	\[
	\begin{tikzcd}
		\Coh[X]\arrow[r]\arrow[d,swap]&\Coh[\bX]\arrow[r]\arrow[d]&\Coh[U]\arrow[d]\\
		\Coh[\X]\arrow[r]&\Coh[\widetilde{\X}]\arrow[r]&\Coh[W]\mathrlap{.}
	\end{tikzcd}
	\]
	Deleting the middle column yields the square of the theorem for $(X,Z)$ and deleting the left column yields the square of the theorem for $(\bX,Z^\prime)$.\qed
\end{corollary}
\begin{remark}
	When $f\colon \bX\to X$ is the blowup of $X$ in $Z$, the closed subscheme $Z\times_X\bX$ is the exeptional divisor of the blowup, which is in particular locally principal in $\bX$. Thus in the diagram of Corollary~\ref{blowupsetup} the right square is now the square of the theorem for a pair with the subscheme locally principal.
	If we can somehow show that the left square is cartesian the proof of the theorem thus reduces to the locally principal case, in which we can try to understand the situation locally by reducing to the affine case with Lemma~\ref{affred}.
	The discussion of this left square takes place in the next section.
\end{remark}
\newpage
\section{The key lemma}
This section deals with the main step towards a proof of the theorem. The present key lemma yields a descent statement in the case where $X$ is affine and $Z$ is cut out by a principal ideal, as well as letting us reduce the general case to the case where $Z$ is locally cut out by such an ideal by showing the left square of Corollary~\ref{blowupsetup} to be cartesian.
For this section fix the following notation:
\begin{itemize}
	\item 
	$X=\Spec{A}$ a Noetherian affine scheme
	\item 
	$Z$ a closed subscheme cut out by $I=(f_1,\dots,f_n)$, with complement $U$
	\item 
	$\comp{A}$ the $I$-adic completion of $A$
	\item 
	$i\colon \Xs\to X$ the canonical map of \emph{affine schemes} induced by $A\to\comp{A}$
	\item 
	$U^\prime$ the pullback $U\times_X\bX$
	\item 
	$\pi\colon \bX\to X$ a morphism of locally Noetherian schemes
	\item 
	$\Es$ the pullback $\bX\times_X\Xs$
\end{itemize}
Thus we have a commutative diagram of the following form in which both squares are cartesian.
\[
\begin{tikzcd}
	\U^\prime\arrow[r,"k^\prime",swap]\arrow[rr,"j^\prime",bend left=30]\arrow[d,"i^{\prime\prime}",swap]&\Es\arrow[r,"\pi^\prime",swap]\arrow[d,"i^\prime"]&\Xs\arrow[d,"i"]\\
	U\arrow[r,"k"]\arrow[rr,"j",bend right=30,swap]&\bX\arrow[r,"\pi"]&X
\end{tikzcd}
\]\vspace*{-3mm}
\begin{lemma}[Key lemma]\label{keylemma}
	Assume that $\pi$ is qcqs and an isomorphism over $U$.
	Then pulling back along $i$ and $\pi$ induces an equivalence
	\[
	\Psi_X\colon\QCoh[X]\to\QCoh[\Xs]\times_{\QCoh[\Es]}\QCoh[\bX]\mathrlap{.}
	\]
\end{lemma}
This is not quite yet cartesianness of the left square of Corollary~\ref{blowupsetup}, but it is a version where we have replaced $\X=\Spf{A}$ by $X^\prime=\Spec{\comp{A}}$.
We collect some general lemmata in preparation for the proof of the key lemma.
\begin{lemma}\label{exactcartiso}
	Let $B$ be a ring and let
	\[
	\begin{tikzcd}
		0\arrow[r]&M_1\arrow[r,"f_M"]\arrow[d,"\alpha_1",swap]&M_2\arrow[r,"g_M"]\arrow[d,"\alpha_2"]&M_3\arrow[r]\arrow[d,"\alpha_3"]&0\\
		0\arrow[r]&N_1\arrow[r,"f_N"]&N_2\arrow[r,"g_N"]&N_3\arrow[r]&0
	\end{tikzcd}
	\]
	be a commutative diagram in an abelian category. Assume that the upper row is a complex and the lower row is exact. Then two of the following three conditions imply the third:
	\begin{itemize}
		\item 
		the upper row is exact as well,
		\item 
		the morphism $\alpha_1$ is an isomorphism,
		\item 
		the right square is cartesian.
	\end{itemize}
\end{lemma}
\begin{proof}
	First let us assume that both rows are exact.
	The kernel of a map $f\colon M\to N$ is the pullback via $0\to N$. Thus we obtain the following diagram
	\[
	\begin{tikzcd}
		M_1\arrow[r]\arrow[d,"f_M",swap]&0\arrow[d]\\
		M_2\arrow[r,"g_M"]\arrow[d,"\alpha_2",swap]&M_3\arrow[d,"\alpha_3"]\\
		N_2\arrow[r,"g_N"]&N_3\mathrlap{,}
	\end{tikzcd}
	\]
	in which the top square is cartesian. Since $g_M$ is epi, the top square is also cocartesian.
	The map $\alpha_1$ between the kernels of $g_M$ and $g_N$ corresponds to the natural morphism from $M_1$ to the pullback of $g_N$ by the zero morphism, which is the pullback of the large outer cospan.
	Now assume $\alpha_1$ is an isomorphism, which equivalently means that the large outer square is cartesian.
	It is cocartesian as well, since $g_N$ is epi. Thus the bottom square is also cocartesian.
	The kernel of $(\alpha_2,g_M)\colon M_2\to N_2\oplus M_3$ factors through $M_1$, so we obtain a commutative diagram:
	
	\[
	\begin{tikzcd}
		&M_1\arrow[r,"\alpha_1"]\arrow[r,"\cong",swap]\arrow[d,hook,"f_M"]&N_1\arrow[dr,hook,"f_N"]\\
		\ker(\alpha_2,g_M)\arrow[r,hook]\arrow[ur,hook]&M_2\arrow[r,"{(\alpha_2,g_M)}"]&N_2\oplus M_3\arrow[r]&N_2
	\end{tikzcd}
	\]
	showing that $\ker(\alpha_2,g_M)$ is zero, and $(\alpha_2,g_M)$ is mono. Thus the bottom square is cartesian.
	
	Assume that the bottom square is cartesian instead. Then the large square becomes cartesian by pasting, which means that $\alpha_1$ is an isomorphism.
	
	Let us finally assume that both the square is cartesian and $\alpha_1$ is an isomorphism.
	Exactness at $M_1$ is okay since $\alpha_1$ is an iso, while exactness at $M_3$ comes from exactness at $N_3$ and the cartesianness of the square.
	By assumption the upper row is a complex, so it is enough to check that $\ker g_M$ factors through $\im f_M$, which is a short diagram chase.
\end{proof}
The proof idea of the following lemma is taken from \cite{BL}.
\begin{lemma}
	Let $M$ be an $A$-module such that $M[f_i^{-1}]=0$ for every $i$. Then the map
	\[
	\id\otimes 1\colon M\to M\otimes_A\comp{A}
	\]
	is an isomorphism of $A$-modules.
\end{lemma}
\begin{proof}
	For every element $x\in M$ there exists an $n$ such that $f_i^nx=0$ in $M$ for every $i$, since $M[f^{-1}]=0$. In particular we can also choose an $m$ such that $I^mx=0$. Thus we can write $M$ as a sequential colimit over the submodules $M_m$ annihilated by $I^m$. Since the tensor product commutes with colimits it is enough to show that the map is an isomorphism for every $M_m$. Thus we can assume that $I^m$ annihilates $M$.
	Now we can write the map $\id\otimes 1$ as the following composition of isomorphisms
	\[
	\begin{aligned}
		M&\cong M\otimes_A A/I^m\\
		&\cong M\otimes_A \comp{A}/I^m\comp{A}\\
		&\cong M/I^m\otimes_A\comp{A}\\
		&= M\otimes_A\comp{A}\mathrlap{,}
	\end{aligned}
	\]
	where we have used Noetherianity of $A$ for the isomorphism $A/I^m\cong\comp{A}/I^m\comp{A}$.
\end{proof}
\begin{corollary}\label{isoglobal}
	Let $M$ be an $\comp{A}$-module such that $M[f_i^{-1}]=0$ for every $i$. Then the map
	\[
	m\colon M\otimes_A\comp{A}\to M
	\]
	is an isomorphism.
\end{corollary}
\begin{proof}
	This is an immediate consequence of the last lemma.
\end{proof}
\begin{remark}\label{iso}
	Corollary~\ref{isoglobal} equivalently shows the following:
	Let $G$ be a quasicoherent sheaf on $X^\prime$ such that $j^*i_*G=0$. Then the counit $i^*i_*G\to G$ of the pullback-pushforward-adjunction for $i$ is an isomorphism.
\end{remark}
\begin{lemma}\label{faithful}
	The morphism
	\[
	(i,j)\colon \Xs\coprod U\to X
	\]
	is an fpqc cover of $X$.
\end{lemma}
\begin{proof}
	The morphism is certainly quasicompact and it is flat by Noetherianity of $A$.
	It is also surjective since $\Xs$ contains $Z$ as a closed subscheme which is mapped isomorphically to $Z$ in $X$. Thus $(i,j)$ is quasicompact and faithfully flat, so it is fpqc.
\end{proof}
Now we have assembled everything we need to prove the key lemma.
\begin{proof}[Proof of lemma 3.1]
	Recall that we want to show that
	\[
	\Psi_X\colon\QCoh[X]\to\QCoh[\Xs]\times_{\QCoh[\Es]}\QCoh[\bX]\mathrlap{,}
	\]
	which is given by pulling back along $i$ and $\pi$ is an equivalence.
	
	First let us check essential surjectivity.
	Let 
	\[
	(F,G,\varphi)\in \QCoh[\Xs]\times_{\QCoh[\Es]}\QCoh[\bX]
	\]
	be an object, i.e.\@ a quasicoherent sheaf $F$ on $\Xs$ and a quasicoherent sheaf $G$ on $\bX$ together with an isomorphism $\varphi\colon \pi^{\prime*}F\to i^{\prime*}G$.
	The units of the adjunctions between the pushforward and the pullback yield morphisms
	\[
	i_*u_F\colon i_*F\to i_* \pi^\prime_* \pi^{\prime*} F
	\]
	and
	\[
	\pi_*u_G\colon\pi_*G\to \pi_*i^\prime_*i^{\prime*}G=i_*\pi^\prime_*i^{\prime*}G
	\]
	in $\QCoh[X]$.
	Composing with $i_*\pi^\prime_*\varphi^{-1}$ one obtains a morphism
	\[
	\alpha=(i_*u_F,i_*\pi^\prime_*\varphi^{-1}\circ\pi_*u_G)\colon i_*F\oplus \pi_*G\to i_*\pi^\prime_*\pi^{\prime*}F
	\]
	between quasicoherent sheaves ($i$, $\pi$ and $\pi^\prime$ are qcqs, so pushforward along these maps preserves quasicoherence).
	Consider the kernel $M$ of $\alpha$, fitting into an exact sequence
	\[
	0\to M\to i_*F\oplus \pi_*G\xrightarrow{\alpha} i_*\pi^\prime_*\pi^{\prime*}F,
	\]
	where $M$ is again a quasicoherent sheaf on $X$.
	This $M$ is our candidate for a quasicoherent sheaf mapping to $(F,G,\varphi)$ up to isomorphism.
	Let us first check that the pullback of $M$ to $\Xs$ is isomorphic to $F$ and the pullback of $M$ to $U$ is isomorphic to $k^*G$.
	
	First consider the pullback of $M$ to $U$.
	The following square is cartesian because $\pi$ is an isomorphism over $U$:
	\[
	\begin{tikzcd}
		U\arrow[r,equal]\arrow[d,"k",swap]&U\arrow[d,"j"]\\
		\bX\arrow[r,"\pi"]&X
	\end{tikzcd}
	\]
	As $j$ is flat applying flat basechange yields that the morphism
	\[
	k^*c_G\colon k^*\pi^*\pi_* G\to k^*G
	\]
	is an isomorphism.
	
	The pullback via $j$ of the exact sequence defining $M$ is 
	\[
	\begin{tikzcd}
		0\arrow[r]&j^*M\arrow[r]& j^*i_*F\oplus j^*\pi_*G\arrow[r,"j^*\alpha"]& j^*i_*\pi^\prime_*\pi^{\prime*}F\mathrlap{.}
	\end{tikzcd}
	\]
	We can compose the flat basechange isomorphism with the projection to obtain a map
	\[
	\begin{tikzcd}
		\beta_U\colon j^*M\arrow[r]& j^*i_*F\oplus j^*\pi_*G\arrow[r,"\pr_2"]& j^*\pi_*G\arrow[r,equal]&k^*\pi^*\pi_*G\arrow[r,"k^*c_G"]\arrow[r,"\cong",swap]& k^*G\mathrlap{,}
	\end{tikzcd}
	\]
	which we will show to be an isomorphism.
	Note that this map is the pullback under $k$ of the map
	\[
	\begin{tikzcd}
		\beta_{\bX}\colon \pi^*M\arrow[r]& \pi^*i_*F\oplus \pi^*\pi_*G\arrow[r,"\pr_2"]& \pi^*\pi_*G\arrow[r,"c_G"]&G\mathrlap{,}
	\end{tikzcd}
	\]
	which will be one half of the eventual isomorphism $\beta\colon\Psi(M)\to(F,G,\varphi)$.
	
	Since $j$ is flat the pullback along $j$ of the exact sequence defining $M$ is still exact.
	Thus the map $\beta_U$ is an isomorphism if
	\[
	j^*\alpha_1=j^*i_*u_F\colon j^*i_*F\to j^*i_*\pi^\prime_*\pi^{\prime*}F
	\]
	is an isomorphism.
	
	Flat basechange for the square
	\[
	\begin{tikzcd}
		U^\prime\arrow[d,"j^\prime",swap]\arrow[r,"i^{\prime\prime}"]&U\arrow[d,"j"]\\
		\Xs\arrow[r,"i"]&X
	\end{tikzcd}
	\]
	yields a natural isomorphism
	\[
	\sigma_F\colon j^*i_*F\to i^{\prime\prime}_*j^{\prime*}F=i^{\prime\prime}_*k^{\prime*}\pi^{\prime*}F
	\]
	and we obtain the following commutative diagram:
	\[
	\begin{tikzcd}[column sep=45]
		j^*i_*F\arrow[r,"j^*i_*u_F"]\arrow[d,"\cong"]\arrow[d,"\sigma_F",swap]&j^*i_*\pi^{\prime}_*\pi^{\prime*}F\arrow[d,"\cong",swap]\arrow[d,"\sigma_{\pi^\prime_*\pi^{\prime*}F}"]\\
		i^{\prime\prime}_*k^{\prime*}\pi^{\prime*}F\arrow[r,"i^{\prime\prime}_*k^{\prime*}\pi^{\prime*}u_F"]&i^{\prime\prime}_*k^{\prime*}\pi^{\prime*}\pi^{\prime}_*\pi^{\prime*}F
	\end{tikzcd}
	\]
	But $\pi^{\prime*}u_F$ is an isomorphism by the triangle identities for the pullback-pushforward adjunction of $\pi^\prime$, so $j^*i_*u_F$ is also an isomorphism. Thus $\beta_U\colon j^*M\to k^*G$ is an isomorphism as well.
	
	Next consider the pullback of the sequence to $\Xs$
	\[
	\begin{tikzcd}
		0\arrow[r]& i^*M\arrow[r]& i^*i_*F\oplus i^*\pi_*G\arrow[r]& i^*i_*\pi^\prime_*\pi^{\prime*}F\mathrlap{,}
	\end{tikzcd}
	\]
	which is still exact by flatness of $i$.
	We once again obtain a map
	\[
	\begin{tikzcd}
		\beta_{\Xs}\colon i^*M\arrow[r]& i^*i_*F\oplus i^*\pi_*G\arrow[r,"\pr_1"]& i^*i_*F\arrow[r,"c_F"]&F\mathrlap{.}
	\end{tikzcd}
	\]
	This will be the other half of the isomorphism $\beta\colon\Psi(M)\to(F,G,\varphi)$.
	Now we will show that $\beta_{\Xs}$ is an isomorphism.
	
	Flat basechange for the square
	\[
	\begin{tikzcd}
		\Es\arrow[r,"\pi^\prime"]\arrow[d,"i^\prime",swap]&\Xs\arrow[d,"i"]\\
		\bX\arrow[r,"\pi"]&X
	\end{tikzcd}
	\]
	yields a natural isomorphism
	\[
	\tau_G\colon i^*\pi_*G\to\pi^\prime_*i^{\prime*}G,
	\]
	which is explicitly given by
	\[
	\begin{tikzcd}[column sep=35]
		i^*\pi_*G\arrow[r,"i^*\pi_*u_G"]&i^*\pi_*i^\prime_*i^{\prime*}G\arrow[r,equal]&
		i^*i_*\pi^\prime_*i^{\prime*}G\arrow[r,"c_{\pi^\prime_*i^{\prime*}G}"]&\pi^\prime_*i^{\prime*}G\mathrlap{.}
	\end{tikzcd}
	\]
	We have an isomorphism
	\[
	\begin{tikzcd}[column sep=35]
		i^*i_*\pi^\prime_*\pi^{\prime*}F\arrow[r,"i^*i_*\pi^\prime_*\varphi"]\arrow[r,"\cong",swap]&i^*i_*\pi^\prime_*i^{\prime*}G&i^*i_*i^*\pi_*G\arrow[l,"i^*i_*\tau_G",swap]\arrow[l,"\cong"]
	\end{tikzcd}
	\]
	and composing the sequence defining $i^*M$ with this yields an exact sequence
	\[
	\begin{tikzcd}
		0\arrow[r]& i^*M\arrow[r]& i^*i_*F\oplus i^*\pi_*G\arrow[r,"\psi"]& i^*i_*i^*\pi_*G\mathrlap{.}
	\end{tikzcd}
	\]
	Here $\psi$ is explicitly given by
	\[
	(i^*i_*(\tau_G^{-1}\circ \pi^\prime_*\varphi\circ u_F){,}\,i^*i_*\tau_G^{-1}\circ i^*\pi_*u_G)\colon i^*i_*F\oplus i^*\pi_*G\to i^*i_*i^*\pi_*G\mathrlap{.}
	\]
	Let $F=\widetilde{N}$ and let $\pi_*G=\widetilde{L}$, where $N$ is an $\comp{A}$-module and $L$ is an $A$-module.
	With this notation evaluating the sequence on global sections yields an exact sequence
	\[
	\begin{tikzcd}
		0\arrow[r]&\Gamma(\Xs,i^*M)\arrow[r]&N\otimes_{A}\comp{A}\oplus L\otimes_{A}\comp{A}\arrow[r,"\psi"]&L\otimes_{A}\comp{A}\otimes_{A}\comp{A}\mathrlap{.}
	\end{tikzcd}
	\]
	Let us analyse the map $\psi$.
	The first component of $\psi$ is $i^*i_*$ applied to a morphism $F\to i^*\pi_*G$.
	This is on global sections given by $-\otimes_A\comp{A}$ applied to some $\comp{A}$-linear $N\to L\otimes_{A}\comp{A}$, which we denote by $f$.
	Thus the first component of $\psi$ is given by
	\[
	f\otimes \id\colon N\otimes_{A}\comp{A}\to L\otimes_{A}\comp{A}\otimes_{A}\comp{A}\mathrlap{.}
	\]
	The second component $\psi_2$ is by definition the pullback along $i$ of the map
	\[
	\iota_L\coloneqq i_*\tau_G^{-1}\circ\pi_*u_G\colon \pi_*G\to i_*i^*\pi_*G\mathrlap{,}
	\]
	so we have
	\[
	\psi_2=\iota_L\otimes\id\colon L\otimes_{A}\comp{A}\to L\otimes_{A}\comp{A}\otimes_{A}\comp{A}
	\]
	with $\iota_L$ $A$-linear.
	Now we want to show that $\psi_2$ is split by the multiplication map
	\[
	m\colon L\otimes_{A}\comp{A}\otimes_{A}\comp{A}\to L\otimes_{A}\comp{A}\mathrlap{.}
	\]
	Consider the following diagram:
	\[
	\begin{tikzcd}
		i^*\pi_*G\arrow[r,"i^*\pi_*u_G"]\arrow[rrd,"\tau_G",swap,bend right=13]&i^*\pi_*i^\prime_*i^{\prime*}G\arrow[r,equal]&i^*i_*\pi^\prime_*i^{\prime*}G\arrow[d,"c_{\pi^\prime_*i^{\prime*}G}"]&i^*i_*i^*\pi_*G\arrow[l,"i^*i_*\tau",swap]\arrow[d,"c_{i^*\pi_*G}"]\arrow[l,"\cong"]\\
		&&\pi^\prime_*i^{\prime*}G&i^*\pi_*G\arrow[l,"\tau_G",swap]\arrow[l,"\cong"]
	\end{tikzcd}
	\]
	The left triangle commutes by the definition of $\tau$ and the square commutes by naturality of the counit.
	On global sections the top row becomes $\psi_2$ and the right vertical map becomes the multiplication map. Thus this multiplication map splits $\psi_2$.
	
	Now let us show that $\psi$ is surjective.
	In the following commutative diagram the map on cokernels is an isomorphism by a short diagram chase:
	\[
	\begin{tikzcd}
		N\otimes_{A}\comp{A}\oplus L\otimes_{A}\comp{A}\arrow[r,"\psi"]\arrow[d,"\pr_2",swap]&L\otimes_{A}\comp{A}\otimes_{A}\comp{A}\arrow[r]\arrow[d,"\pr"]&\coker(\psi)\arrow[d,"\cong",dashed]\arrow[r]&0\\
		L\otimes_{A}\comp{A}\arrow[r,"\pr\circ(\iota_L\otimes\id)"]&\coker(f\otimes\id)\arrow[r]&\coker\arrow[r]&0
	\end{tikzcd}
	\]
	Consider the following diagram:
	\[
	\begin{tikzcd}
		L\otimes_{A}\comp{A}\arrow[r,"\pr\circ(\iota_L\otimes\id)"]\arrow[d,"\iota_L\otimes\id",swap]&\coker(f\otimes\id)\arrow[d,"\cong"]\\
		L\otimes_{A}\comp{A}\otimes_{A}\comp{A}\arrow[r,"\pr\otimes\id",two heads]\arrow[d,"m",swap]&\coker(f)\otimes_{A}\comp{A}\arrow[d,"m"]\\
		L\otimes_{A}\comp{A}\arrow[r,"\pr",two heads]&\coker(f)
	\end{tikzcd}
	\]
	This diagram clearly commutes and the left composition is the identity on $L\otimes_{A}\comp{A}$ by the discussion above.
	If the multiplication map on $\coker(f)$ is an isomorphism as well, the induced map
	\[
	\coker(L\otimes_{A}\comp{A}\xrightarrow{\pr\circ(\iota_L\otimes\id)}\coker(f\otimes\id))\to\coker(L\otimes_{A}\comp{A}\xrightarrow{\pr}\coker(f))=0
	\]
	is an isomorphism, so $\psi$ is surjective.
	
	Thus we are left to show that this multiplication map is an isomorphism.
	For this consider $j^*i_*f$, which fits into the following diagram:
	\[
	\begin{tikzcd}[column sep=50]
		j^*i_*F\arrow[r,"j^*i_*f"]\arrow[d,"\sigma_F",swap]\arrow[d,"\cong"]&j^*i_*i^*\pi_*G\arrow[d,"\sigma_{i^*\pi_*G}"]\arrow[d,"\cong",swap]\\
		i^{\prime\prime}_*k^{\prime*}\pi^{\prime*}F\arrow[r,"i^{\prime\prime}_*k^{\prime*}\pi^{\prime*}f"]\arrow[d,"i^{\prime\prime}_*k^{\prime*}\pi^{\prime*}u_F",swap]&i^{\prime\prime}_*k^{\prime*}\pi^{\prime*}i^*\pi_*G\\
		i^{\prime\prime}_*k^{\prime*}\pi^{\prime*}\pi^\prime_*\pi^{\prime*}F\arrow[r,"i^{\prime\prime}_*k^{\prime*}\pi^{\prime*}\pi^\prime_*\varphi"]\arrow[r,swap,"\cong"]&i^{\prime\prime}_*k^{\prime*}\pi^{\prime*}\pi^\prime_*i^{\prime*}G\arrow[u,"i^{\prime\prime}_*\pi^{\prime*}k^{\prime*}\tau_G^{-1}",swap]\arrow[u,"\cong"]
	\end{tikzcd}
	\]
	The top square commutes by naturality of $\sigma$ and the bottom square commutes by the definition of $f$.
	By the triangle identities $\pi^{\prime*}u_F$ is an isomorphism, so $j^*i_*f$ is an isomorphism as well.
	Thus 
	\[
	m\colon\coker(f)\otimes_A\comp{A}\to\coker(f)
	\]
	is an isomorphism by Remark~\ref{iso} and $\psi$ is surjective.
	
	Denote the image of $f$ by 
	\[
	\nu\colon N^\prime\hookrightarrow L\otimes_{A}\comp{A}
	\]
	and the kernel of $f$ by
	\[
	\kappa\colon K\to N\mathrlap{.}
	\]
	By Noetherianness of $A$ the completion $\comp{A}$ is flat over $A$, so
	\[
	\nu\otimes\id\colon N^\prime\otimes_{A}\comp{A}\hookrightarrow L\otimes_{A}\comp{A}\otimes_{A}\comp{A}
	\]
	is injective as well. Denote the map
	\[
	(\nu\otimes\id\,{,}\,\iota_L\otimes\id)\colon N^\prime\otimes_{A}\comp{A}\oplus L\otimes_{A}\comp{A}\to L\otimes_{A}\comp{A}\otimes_{A}\comp{A}
	\]
	by $\psi^\prime$. This fits into a commutative diagram
	\[
	\begin{tikzcd}
		N\otimes_{A}\comp{A}\oplus L\otimes_{A}\comp{A}\arrow[d,"(f\otimes\id\,{,}\,\id)",swap]\arrow[r,"\psi"]&L\otimes_{A}\comp{A}\otimes_{A}\comp{A}\arrow[d,equal]\\
		N^{\prime}\otimes_{A}\comp{A}\oplus L\otimes_{A}\comp{A}\arrow[r,"\psi^\prime"]&L\otimes_{A}\comp{A}\otimes_{A}\comp{A}\mathrlap{.}
	\end{tikzcd}
	\]
	Consider the following commutative diagram with exact rows:
	\[
	\begin{tikzcd}
		0\arrow[r]&\ker(\psi^\prime)\arrow[r]\arrow[d,dashed,"\cong",swap]&N^\prime\otimes_{A}\comp{A}\oplus L\otimes_{A}\comp{A}\arrow[r,"\psi^\prime"]\arrow[d,"\pr_2"]& L\otimes_{A}\comp{A}\otimes_{A}\comp{A}\arrow[r]\arrow[d,"\pr"]&0\\
		0\arrow[r]&\ker\arrow[r]&L\otimes_{A}\comp{A}\arrow[r,"\pr\circ(\iota_L\otimes\id)"]& \coker(\nu\otimes\id)\arrow[r]&0
	\end{tikzcd}
	\]
	Using exactness and the injectivity of $\nu\otimes\id$ a diagram chase shows that the dashed induced map on kernels is an isomorphism.
	A similar argument to the argument for the cokernel earlier yields an isomorphism over $L\otimes_A\comp{A}$:
	\[
	\ker(L\otimes_A\comp{A}\xrightarrow{\pr\circ(\iota_L\otimes\id)}
	\coker(\nu\otimes\id))\cong\ker(L\otimes_A\comp{A}\xrightarrow{\pr}\coker(\nu))
	\]
	Since $\nu$ is injective this kernel is isomorphic to $\nu$ itself.
	Plugging this into the exact sequence of $\psi^\prime$ yields an exact sequence
	\[
	\begin{tikzcd}
		0\arrow[r]
		&N^\prime\arrow[r,"(\iota_{N^\prime}{,}\,\nu)"]
		&N^\prime\otimes_{A}\comp{A}\oplus L\otimes_{A}\comp{A}\arrow[r,"\psi^\prime"]
		&L\otimes_{A}\comp{A}\otimes_{A}\comp{A}\arrow[r]&0
	\end{tikzcd}
	\]
	for some $A$-linear map
	\[
	\iota_{N^\prime}\colon N^\prime\to N^\prime\otimes_A\comp{A}\mathrlap{.}
	\]
	The following diagram commutes since $(\iota_{N^\prime},\,\nu)$ is the kernel of $\psi^\prime$
	\[
	\begin{tikzcd}
		N^\prime\arrow[r,"\iota_{N^\prime}"]\arrow[d,"\nu",swap,hook]&N^\prime\otimes_{A}\comp{A}\arrow[d,"\nu\otimes\id",hook]\arrow[r,"m"]&N^\prime\arrow[d,"\nu",hook]\\
		L\otimes_{A}\comp{A}\arrow[r,"\iota_L\otimes\id"]\arrow[rr,"\id",bend right=20,swap]&L\otimes_{A}\comp{A}\otimes_{A}\comp{A}\arrow[r,"m"]&L\otimes_A\comp{A}\mathrlap{,}
	\end{tikzcd}
	\]
	and as $\iota_L\otimes\id$ is split by $m$ the bottom map is the identity. Thus the map $\iota_{N^\prime}$ is split by the multiplication map on $N^\prime$.
	
	The sequence above fits into a commutative diagram
	\[
	\begin{tikzcd}
		0\arrow[r]&\Gamma(\Xs,i^*M)\arrow[r]\arrow[d,dashed]&N\otimes_{A}\comp{A}\oplus L\otimes_{A}\comp{A}\arrow[r,"\psi"]\arrow[d,"(f\otimes\id\,{,}\,\id)"]& L\otimes_{A}\comp{A}\otimes_{A}\comp{A}\arrow[r]\arrow[d,equal]&0\\
		0\arrow[r]&N^\prime\arrow[r,"(\iota_{N^\prime}{,}\,\nu)"]&N^\prime\otimes_{A}\comp{A}\oplus L\otimes_{A}\comp{A}\arrow[r,"\psi^\prime"]& L\otimes_{A}\comp{A}\otimes_{A}\comp{A}\arrow[r]&0\mathrlap{.}
	\end{tikzcd}
	\]
	The rightmost vertical morphism is an isomorphism and the rows are exact, so the left square is cocartesian by Lemma~\ref{exactcartiso}. Since $\Gamma(\Xs,i^*M)\to N\otimes_{A}\comp{A}\oplus L\otimes_{A}\comp{A}$ is mono the square is cartesian as well. Thus the following square is also cartesian:
	\[
	\begin{tikzcd}
		\Gamma(\Xs,i^*M)\arrow[r]\arrow[d]&N\otimes_{A}\comp{A}\arrow[d,"f\otimes\id"]\\
		N^\prime\arrow[r,"\iota_{N^\prime}"]&N^\prime\otimes_{A}\comp{A}
	\end{tikzcd}
	\]
	We have seen earlier that $j^*i_*f$ is an isomorphism, so $j^*i_*\ker(f)$ is zero. Thus
	\[
	m\colon K\otimes_{A}\comp{A}\to K
	\]
	is an isomorphism by Remark~\ref{iso}.
	Using this map we can consider the following diagram:
	\[
	\begin{tikzcd}
		0\arrow[r]&K\arrow[r,dashed]\arrow[d,"m^{-1}",swap]\arrow[d,"\cong"]&\Gamma(\Xs,i^*M)\arrow[r]\arrow[d]\arrow[dr, phantom, "\scalebox{1.3}{$\lrcorner$}", very near start, xshift=-9pt, color=black]&N^\prime\arrow[r]\arrow[d,"\iota_{N^\prime}"]&0\\
		0\arrow[r]&K\otimes_{A}\comp{A}\arrow[r,"\kappa\otimes\id"]\arrow[d,"m",swap]\arrow[d,"\cong"]&N\otimes_{A}\comp{A}\arrow[r,"f\otimes\id"]\arrow[d,"m"]&N^\prime\otimes_{A}\comp{A}\arrow[r]\arrow[d,"m"]&0\\
		0\arrow[r]&K\arrow[r,"\kappa"]&N\arrow[r,"f"]&N^\prime\arrow[r]&0
	\end{tikzcd}
	\]
	
	The dashed map is induced by $((\kappa\otimes\id)\circ m^{-1},0)$ via the universal property of the pullback, making the diagram commutative.
	The top row is certainly a complex, $m^{-1}$ is an isomorphism and the top right square is cartesian so we can use Lemma~\ref{exactcartiso} to conclude that the top row is exact as well.
	
	Both compositions $K\to K$ and $N^\prime\to N^\prime$ are the identity, so the 5-lemma shows that the central composition $\Gamma(\Xs,i^*M)\to N$ is an isomorphism.
	When retracing the definitions this map is the map induced on global sections by
	\[
	\begin{tikzcd}
		\beta_{\Xs}\colon i^*M\arrow[r]&i^*i_*F\oplus i^*\pi_*G\arrow[r,"\pr_1"]&i^*i_*F\arrow[r,"c_F"]&F\mathrlap{,}
	\end{tikzcd}
	\]
	so $\beta_{\Xs}$ is an isomorphism.
	
	As a next step we want to show that $(\beta_{\Xs},\beta_{\bX})$ is a morphism $\Psi(M)\to(F,G,\varphi)$ in the pullback category.
	For this we need to check that the pullbacks of $\beta_{\bX}$ and $\beta_{\Xs}$ to $\Es$ make the following square commute:
	\[
	\begin{tikzcd}[column sep=40]
		\pi^{\prime*}i^*M\arrow[r,"\pi^{\prime*}\beta_{\Xs}"]\arrow[d,"\cong",swap]&\pi^{\prime*}F\arrow[d,"\varphi"]\\
		i^{\prime*}\pi^*M\arrow[r,"i^{\prime*}\beta_{\bX}"]&i^{\prime*}G
	\end{tikzcd}
	\]
	First consider the following diagram:
	\[
	\begin{tikzcd}[column sep=40]
		\pi^{\prime*}i^*\pi_*G\arrow[r,"\pi^{\prime*}i^*\pi_*u_G"]\arrow[d,"\cong",swap]&
		\pi^{\prime*}i^*\pi_*i^\prime_*i^{\prime*}G\arrow[r,equal]\arrow[d,"\cong"]&
		\pi^{\prime*}i^*i_*\pi^\prime_*i^{\prime*}G\arrow[r,"\pi^{\prime*}c_{\pi^\prime_*i^{\prime*}G}"]&
		\pi^{\prime*}\pi^\prime_*i^{\prime*}G\arrow[d,"c_{i^{\prime*}G}"]\\
		i^{\prime*}\pi^*\pi_*G\arrow[r,"i^{\prime*}\pi^*\pi_*u_G"]\arrow[drr,"i^{\prime*}c_G",swap]&
		i^{\prime*}\pi^*\pi_*i^\prime_*i^{\prime*}G\arrow[r,"i^{\prime*}c_{i^\prime_*i^{\prime*}G}"]&
		i^{\prime*}i^\prime_*i^{\prime*}G\arrow[r,"c_{i^{\prime*}G}"]&
		i^{\prime*}G\\
		&&i^{\prime*}G\arrow[u,"i^{\prime*}u_G"]\arrow[ur,"\id_{i^{\prime*}G}",swap,bend right=18]
	\end{tikzcd}
	\]
	Here the top composition is $\pi^{\prime*}\tau_G$ and everything but the large square clearly commutes.
	The large square commutes as the top composition
	\[
	\pi^{\prime*}i^*i_*\pi^\prime_*i^{\prime*}G\to i^{\prime*}G
	\]
	is the counit of $i^{\prime*}G$ of the adjunction $\pi^{\prime*}i^*\dashv i_*\pi^\prime_*$, while the bottom composition
	\[
	i^{\prime*}\pi^*\pi_*i^\prime_*i^{\prime*}G\to i^{\prime*}G
	\]
	is the counit of $i^{\prime*}G$ of the adjunction $i^{\prime*}\pi^*\dashv\pi_*i^\prime_*$.
	There are isomorphisms $\pi^{\prime*}i^*\cong i^{\prime*}\pi^*$ and $i_*\pi^\prime_*\cong\pi_*i^\prime_*$, so the counits agree after applying these isomorphisms, so the square commutes.
	Thus the diagram is commutative.	
	Now we can consider 
	\[
	\begin{tikzcd}[column sep=16]
		\pi^{\prime*}i^*M\arrow[r]\arrow[d,"\cong",swap]&[-3pt]\pi^{\prime*}i^*i_*F\oplus \pi^{\prime*}i^*\pi_*G\arrow[r,"\pr_2"]\arrow[d,"\cong"]&\pi^{\prime*}i^*\pi_*G\arrow[r,"\pi^{\prime*}\tau_G"]\arrow[d,"\cong"]&[10pt]\pi^{\prime*}\pi^\prime_*i^{\prime*}G\arrow[r,"\pi^{\prime*}\pi^{\prime}_*\varphi^{-1}"]\arrow[d,"c_{i^{\prime*}G}"]&[14pt]\pi^{\prime*}\pi^{\prime}_*\pi^{\prime*}F\arrow[d,"c_{\pi^{\prime*}F}"]\\
		i^{\prime*}\pi^*M\arrow[r]&[-3pt]i^{\prime*}\pi^*i_*F\oplus i^{\prime*}\pi^*\pi_*G\arrow[r,"\pr_2"]&i^{\prime*}\pi^*\pi_*G\arrow[r,"i^{\prime*}c_G"]&[10pt]i^{\prime*}G&[14pt]\pi^{\prime*}F\mathrlap{.}\arrow[l,"\varphi",swap]
	\end{tikzcd}
	\]
	In this diagram everything but the third square clearly commutes, but this third square is the outer square of the diagram that we just looked at, so this diagram is also commutative.	
	The lower composition $i^{\prime*}\pi^*M\to i^{\prime*}G$ is $i^{\prime*}\beta_{\bX}$, so checking that the top composition $\pi^{\prime*}i^*M\to\pi^{\prime*}F$ is $\pi^{\prime*}\beta_{\Xs}$ yields the claim.
	
	For this consider the following diagram:
	\[
	\begin{tikzcd}[column sep=40]
		\pi^{\prime*}i^*M\arrow[r]&[-20pt]
		\pi^{\prime*}i^*i_*F\oplus \pi^{\prime*}i^*\pi_*G\arrow[r,"\pr_2"]\arrow[d,"\pr_1",swap]&
		\pi^{\prime*}i^*\pi_*G\arrow[d,"\pi^{\prime*}i^*\alpha_2"]\arrow[r,"\pi^{\prime*}\tau_G"]&
		\pi^{\prime*}\pi^\prime_*i^{\prime*}G\arrow[d,"\pi^{\prime*}\pi^{\prime}_*\varphi^{-1}"]\\
		&\pi^{\prime*}i^*i_*F\arrow[r,"\pi^{\prime*}i^*i_*u_F"]\arrow[d,"\pi^{\prime*}c_F",swap]&
		\pi^{\prime*}i^*i_*\pi^\prime_*\pi^{\prime*}F\arrow[r,"\pi^{\prime*}c_{\pi^\prime_*\pi^{\prime*}F}"]\arrow[d,"\pi^{\prime*}c_{\pi^\prime_*\pi^{\prime*}F}"]
		&\pi^{\prime*}\pi^\prime_*\pi^{\prime*}F\arrow[d,"c_{\pi^{\prime*}F}"]\\
		&\pi^{\prime*}F\arrow[r,"\pi^{\prime*}u_F"]&
		\pi^{\prime*}\pi^\prime_*\pi^{\prime*}F\arrow[r,"c_{\pi^{\prime*}F}"]&
		\pi^{\prime*}F
	\end{tikzcd}
	\]
	The lower composition from $\pi^{\prime*}i^*M$ to $\pi^{\prime*}F$ is $\pi^{\prime*}\beta_{\Xs}$ and the upper composition is the upper composition from the last diagram. Thus it is enough to check that this diagram commutes.
	The two bottom squares clearly commute, and the top left square commutes after composing with the map from $\pi^{\prime*}i^*M$ since $\pi^{\prime*}i^*i_*u_F$ is $\pi^{\prime*}i^*\alpha_1$ and $M$ is the kernel of $\alpha$.
	For commutativity of the top right square consider the commutative diagram
	\[
	\begin{tikzcd}[column sep=large]
		\pi^{\prime*}i^*\pi_*G\arrow[r,"\pi^{\prime*}i^*\pi_*u_G"]
		&[5pt]\pi^{\prime*}i^*\pi_*i^\prime_*i^{\prime*}G\arrow[r,equal]
		&[-10pt]\pi^{\prime*}i^*i_*\pi^\prime_*i^{\prime*}G\arrow[r,"\pi^{\prime*}c_{\pi^\prime_*i^{\prime*}G}"]\arrow[d,"\pi^{\prime*}i^*i_*\pi^{\prime*}\varphi^{-1}",swap]
		&\pi^{\prime*}\pi^\prime_*i^{\prime*}G\arrow[d,"\pi^{\prime*}\pi^\prime_*\varphi^{-1}"]\\
		&&\pi^{\prime*}i^*i_*\pi^{\prime}_*\pi^{\prime*}F\arrow[r,"\pi^{\prime*}c_{\pi^\prime_*\pi^{\prime*}F}"]
		&\pi^{\prime*}\pi^\prime_*\pi^{\prime*}F\mathrlap{.}
	\end{tikzcd}
	\]
	The composition 
	\[
	\pi^{\prime*}i^*\pi_*G\to\pi^{\prime*}\pi^\prime_*i^{\prime*}G
	\]
	in the diagram is $\pi^{\prime*}\tau_G$ and the composition
	\[
	\pi^{\prime*}i^*\pi_*G\to\pi^{\prime*}i^*i_*\pi^{\prime}_*\pi^{\prime*}F
	\]
	is $\pi^{\prime*}i^*\alpha_2$.
	Thus the large diagram commutes, showing that
	\[
	\begin{tikzcd}[column sep=40]
		\pi^{\prime*}i^*M\arrow[r,"\pi^{\prime*}\beta_{\Xs}"]\arrow[d,"\cong",swap]&\pi^{\prime*}F\arrow[d,"\varphi"]\arrow[d,swap,"\cong"]\\
		i^{\prime*}\pi^*M\arrow[r,"i^{\prime*}\beta_{\bX}"]&i^{\prime*}G
	\end{tikzcd}
	\]
	commutes, so $(\beta_{\Xs},\beta_{\bX})$ is a morphism $\Psi(M)\to(F,G,\varphi)$ in the pullback category.
	
	Now we are left to show that $(\beta_{\Xs},\beta_{\bX})$ is an isomorphism, for which it is enough to check that $\beta_{\Xs}$ and $\beta_{\bX}$ are isomorphisms. We have already seen that $\beta_{\Xs}$ is an isomorphism, so consider $\beta_{\bX}$.
	The pullback of $\beta_{\bX}$ to $U$ is $\beta_U$, which is an isomorphism. By the commutative square above $i^{\prime*}\beta_{\bX}$ is an isomorphism iff $\pi^{\prime*}\beta_{\Xs}$ is an isomorphism. But $\beta_{\Xs}$ is an isomorphism, so the pullback of $\beta_{\bX}$ to $\Es$ is an isomorphism as well.
	By Lemma~\ref{faithful} $(\Xs,U)$ is an fpqc cover of $X$ and $(\Es,U)$ is an fpqc cover of $\bX$ as the pullback of this cover, so $\beta_{\bX}$ is an isomorphism too.
	Thus $(\beta_{\Xs},\beta_{\bX})$ is an isomorphism $\Psi(M)\to(F,G,\varphi)$ in the pullback category and $\Psi$ is essentially surjective.
	
	Let us now show faithfulness and fullness.
	The following diagram is commutative up to natural isomorphism (as $\Psi=(i^*,\pi^*)$):
	\[
	\begin{tikzcd}
		\QCoh[X]\arrow[r,"\Psi"]\arrow[dr,"(i^*{,}\,j^*)",swap]&\QCoh[\Xs]\times_{\QCoh[\Es]}\QCoh[\bX]\arrow[d,"(\id{,}\,k^*)"]\\
		&\QCoh[\Xs]\times\QCoh[U]
	\end{tikzcd}
	\]
	By Lemma~\ref{faithful} $(\Xs,U)$ is an fpqc cover of $X$, so the diagonal functor is faithful. Thus $\Psi$ is faithful as well.
	
	To show fullness of $\Psi$ we will first construct an isomorphism $\eta$ between $M$ and the kernel of the exact sequence 
	\[
	\begin{tikzcd}
		0\arrow[r]&\ker\alpha\arrow[r]&i_*i^*M\oplus \pi_*\pi^*M\arrow[r,"\alpha"]&i_*\pi^{\prime}_*\pi^{\prime*}i^*M
	\end{tikzcd}
	\]
	occuring in the discussion of essential surjectivity for $(F,G,\varphi)=\Psi(M)$.
	Once we have this isomorphism we can construct a morphism between the kernels from a morphism $\Psi(M)\to \Psi(N)$ and compose with $\eta$ to obtain a morphism $M\to N$.
	
	The map $\alpha$ is given by $\alpha_1=i_*u_{i^*M}$ and
	\[
	\begin{tikzcd}
		\alpha_2\colon \pi_*\pi^*M\arrow[r,"\pi_*u_{\pi^*M}"]&\pi_*i^\prime_*i^{\prime*}\pi^*M\arrow[r,"\cong"]&i_*\pi^{\prime}_*\pi^{\prime*}i^*M\mathrlap{.}
	\end{tikzcd}
	\]
	We can define a map 
	\[
	\eta=(-u^i_M,u^\pi_M)\colon M\to i_*i^*M\oplus \pi_*\pi^*M
	\]
	using the units of the adjunctions.
	The composition $-\alpha_1\circ\eta_1$ is
	\[
	i_*u_{i^*M}\circ u^i_M\colon M\to i_*\pi^{\prime}_*\pi^{\prime*}i^*M
	\]
	and $\alpha_2\circ\eta_2$ is
	\[
	\pi_*u_{\pi^*M}\circ u^\pi_M\colon M\to \pi_*i^\prime_*i^{\prime*}\pi^*M
	\]
	up to the isomorphism
	\[
	\pi_*i^\prime_*i^{\prime*}\pi^*M\to i_*\pi^{\prime}_*\pi^{\prime*}i^*M\mathrlap{.}
	\]
	These compositions are the units of the adjunctions $i^{\prime*}\pi^*\dashv\pi_*i^\prime_*$ and $\pi^{\prime*}i^*\dashv i_*\pi^{\prime}_*$, respectively.
	As these functors are isomorphic their units are equal after applying the isomorphism between them, which is exactly the isomorphism in the definition of $\alpha_2$ above.
	Thus $-\alpha_1\circ\eta_1=\alpha_2\circ\eta_2$, so $\eta$ factors through $\ker\alpha$, say as $\eta\colon M\to\ker\alpha$.
	After pulling back to $\Xs$ via $i$ we obtain a diagram
	\[
	\begin{tikzcd}
		&i^*M\arrow[d,"-i^*u_M"]\arrow[dl,"i^*\eta",swap]\\
		i^*\ker\alpha\arrow[r]\arrow[dr,"\beta_{\Xs}",swap]&i^*i_*i^*M\arrow[d,"c_{i^*M}"]\\
		&i^*M\mathrlap{,}
	\end{tikzcd}
	\]
	where both triangles commute by construction and the right composition is $-\id_{i^*M}$. Since $\beta_{\Xs}$ is an isomorphism by the proof of essential surjectivity, so is $i^*\eta$.
	Pulling back to $\bX$ and running the same argument shows that $\pi^*\eta$ is an isomorphism as well, so in particular $j^*\eta$ is also an isomorphism.
	As $(\Xs,U)$ is an fpqc cover of $X$ by Lemma~\ref{faithful}, $\eta$ is an isomorphism.
	
	Now let $M,N$ be quasicoherent sheaves on $X$ and let 
	\[
	(f,g)\colon\Psi(M)\to \Psi(N)
	\]
	be a morphism.
	As mentioned earlier we will now construct a map $M\to N$ by considering a map between $\ker\alpha_M$ and $\ker\alpha_N$ and composing with $\eta_M$ and $\eta_N^{-1}$.
	The morphism $(f,g)$ is explicitly given by maps $f\colon i^*M\to i^*N$ and $g\colon \pi^*M\to\pi^*N$ such that
	\[
	\begin{tikzcd}
		\pi^{\prime*}i^*M\arrow[r,"\pi^{\prime*}f"]\arrow[d,"\cong",swap]&\pi^{\prime*}i^*N\arrow[d,"\cong"]\\
		i^{\prime*}\pi^*M\arrow[r,"i^{\prime*}g"]&i^{\prime*}\pi^*N
	\end{tikzcd}
	\]
	commutes.
	Consider
	\[
	\begin{tikzcd}
		i_*i^*M\oplus\pi_*\pi^*M\arrow[r,"\alpha_M"]\arrow[d,"(i_*f{,}\pi_*g)",swap]&i_*\pi^{\prime}_*\pi^{\prime*}i^*M\arrow[d,"i_*\pi^\prime_*\pi^{\prime*}f"]\\
		i_*i^*N\oplus\pi_*\pi^*N\arrow[r,"\alpha_N"]&i_*\pi^{\prime}_*\pi^{\prime*}i^*N\mathrlap{,}
	\end{tikzcd}
	\]
	where $\alpha_M$ is again given by
	\[
	\begin{tikzcd}[column sep=70]
		i_*i^*M\oplus\pi_*\pi^*M\arrow[r,"(i_*u_{i^*M}{,}\,\pi_*u_{\pi^*M})"]&
		i_*\pi_*\pi^{\prime*}i^*M\oplus\pi_*i^\prime_*i^{\prime*}\pi^*M\arrow[r,"(\id\,{,}\,\cong)"]&[-30pt]
		i_*\pi^{\prime}_*\pi^{\prime*}i^*M.
	\end{tikzcd}
	\]
	When restricting to the first factor the square commutes by naturality of the unit.
	Restricting instead to the second factor it becomes the outer square in the following diagram:
	\[
	\begin{tikzcd}[column sep=35]
		\pi_*\pi^*M\arrow[r,"\pi_*u_{\pi^\prime M}"]\arrow[d,"\pi_*g",swap]&\pi_*i^\prime_*i^{\prime*}\pi^\prime M\arrow[r,"\cong"]\arrow[d,"\pi_*i^\prime_*i^{\prime*}g"]&\pi_*i^\prime_*\pi^{\prime*}i^* M\arrow[r,equal]\arrow[d,"\pi_*i^\prime_*\pi^{\prime*}f"]&i_*\pi^{\prime}_*\pi^{\prime*}i^* M\arrow[d,"i_*\pi^\prime_*\pi^{\prime*}f"]\\
		\pi_*\pi^*M\arrow[r,"\pi_*u_{\pi^\prime N}"]&\pi_*i^\prime_*i^{\prime*}\pi^\prime N\arrow[r,"\cong"]&\pi_*i^\prime_*\pi^{\prime*}i^* N\arrow[r,equal]&i_*\pi^{\prime}_*\pi^{\prime*}i^* N
	\end{tikzcd}
	\]
	The left square commutes by the naturality of the unit and the center square is the pushforward of the commutative square assumed to commute for $f$ and $g$.
	
	Thus we obtain an induced map $\ker\alpha_M\to\ker\alpha_N$, compose this with the isomorphisms $\eta_M$ and $\eta_N^{-1}$ to obtain a map $h\colon M\to N$.
	Now we need to check that the pullback of $h$ to $\Xs$ is $f$ and the pullback of $h$ to $\bX$ is $g$.
	Let us check this for $\bX$, the other calculation is similar.
	Consider the following diagram, which is commutative by the construction of $\eta$ and $h$:
	\[
	\begin{tikzcd}
		\pi^*M\arrow[d,"\pi^*\eta_M"]\arrow[ddd,bend right=50,"\pi^*h",swap]\arrow[dr,"(\pi^*u_M{,}\,\pi^*u_M)"]\arrow[drr,bend left=20,"\pi^*u_M"]\\
		\pi^*\ker\alpha_M\arrow[d]\arrow[r]&\pi^*i_*i^*M\oplus \pi^*\pi_*\pi^*M\arrow[d,"(\pi^*i_*f{,}\,\pi^*\pi_*g)"]\arrow[r,"\pr_2"]&\pi^*\pi_*\pi^*M\arrow[d,"\pi^*\pi_*g"]\arrow[r,"c_{\pi^*M}"]&\pi^*M\arrow[d,"g"]\\
		\pi^*\ker\alpha_N\arrow[r]\arrow[d,"\pi^*\eta_N^{-1}"]&\pi^*i_*i^*N\oplus \pi^*\pi_*\pi^*N\arrow[r,"\pr_2"]&\pi^*\pi_*\pi^*N\arrow[r,"c_{\pi^*N}"]&\pi^*N\\
		\pi^*N\arrow[ur,"(\pi^*u_N{,}\,\pi^*u_N)",swap]\arrow[urr,bend right=20,"\pi^*u_N",swap]
	\end{tikzcd}
	\]
	The triangle identities yield that the top and bottom compositions $\pi^*M\to\pi^*M$ and $\pi^*N\to\pi^*N$ are the identity. Thus $\pi^*h=g$, showing that $\Psi$ is full.
\end{proof}
\begin{remark}
	In the case where $X$ is affine, $Z$ is principal and $\pi\colon U\to X$ is the inclusion the lemma shows that the square
	\[
	\begin{tikzcd}
		\Mod[A]\arrow[r]\arrow[d]&\Mod[\comp{A}]\arrow[d]\\
		\Mod[A{[f^{-1}]}]\arrow[r]&\Mod[\comp{A}{[f^{-1}]}]
	\end{tikzcd}
	\]
	is cartesian, so it recovers a variant of the original formulation of the Beauville-Laszlo theorem in which we have assumed that $A$ is Noetherian, but dropped the regularity assumption on the modules.
\end{remark}
To apply the lemma during the proof of the main theorem we will need a version dealing with coherent sheaves. We can immediately show a more general statement for properties which have fpqc descent.
\begin{corollary}\label{keylemmafaithful}
	Let $\mathcal{P}$ be a property of coherent sheaves on locally Noetherian schemes such that
	\begin{itemize}
		\item 
		for every fpqc cover $\{f_i\colon Y_i\to Z\}$ of locally Noetherian schemes a coherent sheaf $F$ on $Z$ has $\mathcal{P}$ if every $f_i^*F$ has $\mathcal{P}$,
		\item 
		the property $\mathcal{P}$ is stable under pullback.
	\end{itemize}
	Denote the full subcategory of $\QCoh[-]$ of coherent sheaves that satisfy $\mathcal{P}$ by $\QCohp[-]$.
	Under the assumptions of the Key Lemma~\ref{keylemma} pullback via $i$ and $\pi$ also yields an equivalence
	\[
	\Psi_X^{\mathcal{P}}\colon\QCohp[X]\to\QCohp[\Xs]\times_{\QCohp[\Es]}\QCohp[\bX]\mathrlap{.}
	\]
\end{corollary}
\begin{proof}
	The functor $\Psi_X^{\mathcal{P}}$ is the restriction of the $\Psi_X$ of the lemma to $\QCohp[X]$, so we obtain a commutative square
	\[
	\begin{tikzcd}
		&\QCohp[X]\arrow[r,"\Psi_X^{\mathcal{P}}"]\arrow[d,hook]&\QCohp[\Xs]\times_{\QCohp[\Es]}\QCohp[\bX]\arrow[d,hook,shorten >=-.4em]\\
		&\QCoh[X]\arrow[r,"\Psi_X"]\arrow[r,"\cong",swap]&\QCoh[\Xs]\times_{\QCoh[\Es]}\QCoh[\bX]
	\end{tikzcd}
	\]
	in which the lower functor is an equivalence by the key lemma.
	The restriction of $\Psi_X$ to the full subcategory $\QCohp[-]$ is again fully faithful.
	For essential surjectivity let $(F,G,\varphi)$ be an object in $\QCohp[\Xs]\times_{\QCohp[\Es]}\QCohp[\bX]$; we can find an quasicoherent $M$ and an isomorphism
	\[
	\Psi_X(M)\to(F,G,\varphi)
	\] 
	in $\QCoh[\Xs]\times_{\QCoh[\Es]}\QCoh[\bX]$.
	In particular $i^*M$ and $\pi^*M$ are isomorphic to $F$ and $G$, so they have $\mathcal{P}$. Thus $j^*M$ has $\mathcal{P}$, too.
	By Lemma~\ref{faithful} the map $(i,j)\colon\Xs\coprod U\to X$ is an fpqc cover of $X$, so $M$ also has $\mathcal{P}$, showing essential surjectivity.
\end{proof}
\begin{corollary}\label{keylemmacoh}
	Under the assumptions of the Key Lemma~\ref{keylemma} pullback via $i$ and $\pi$ yields an equivalence
	\[
	\Psi^{\operatorname{coh}}_X\colon\Coh[X]\to\Coh[\Xs]\times_{\Coh[\Es]}\Coh[\bX]\mathrlap{.}
	\]
\end{corollary}
\begin{proof}
	This is Corollary~\ref{keylemmafaithful} applied to the property of being of finite type.
\end{proof}
\newpage
\section{Proof of the main theorem}
Now that we have proven the key lemma, which - with some adaptation - can both be used to reduce to the locally principal case and is a main ingredient in the locally principal case we are ready to tackle the proof. Let us first show the locally principal case, and then consider the reduction step afterwards.
\subsection{The locally principal case}
\begin{lemma}\label{principal}
	Let $X$ be a locally Noetherian scheme, let $Z$ be a closed subscheme of $X$ which is locally given by the vanishing of a principal ideal.
	
	Then the theorem holds for $(X,Z)$.
\end{lemma}
\begin{proof}
	Let us first assume that $X$ is affine and $Z$ is principal.
	In this case the statement reduces to Corollary~\ref{keylemmacoh} after a few observations.
	Recalling the setup of the theorem we need to show that the outer square of
	\[
	\begin{tikzcd}
		&\Coh[X]\arrow[r,"j^*"]\arrow[d,"s_X^*"]\arrow[ddl,"i^*",bend right=20,swap]&\Coh[U]\arrow[d,"(s_X|_{\U})^*"]\\
		&\Coh[X^{\ad}]\arrow[r]\arrow[d,"i^{\ad,*}"]&\Coh[\U]\arrow[d]\\
		\Coh[\X]\arrow[r,"p^*"]&\Coh[\X^{\ad}]\arrow[r]&\Coh[W]\mathrlap{.}
	\end{tikzcd}
	\]
	is cartesian.
	The functor $p^*$ is an equivalence by Lemma~\ref{pequiv}, $s_X^*$ is an equivalence by Corollary~\ref{sequiv}, moreover $(s_X|_{\U})^*$ is an equivalence by combining Corollary~\ref{sequiv} and Proposition~\ref{isoU}.
	Thus we can alternatively check that the square
	\[
	\begin{tikzcd}
		\Coh[X^{\ad}]\arrow[r]\arrow[d,"i^{\ad,*}",swap]&\Coh[\U]\arrow[d]\\
		\Coh[\X^{\ad}]\arrow[r]&\Coh[W]
	\end{tikzcd}
	\]
	is cartesian.
	But we have already calculated the spaces $\U$ and $W$ in Example~\ref{distopen} and Example~\ref{distopenII}. We have
	\[
	\U=\Spa{A[f^{-1}]}{A^{\cl}}=\R_{X^{\ad}}\big(\frac{f}{f}\big)
	\]
	and
	\[
	W=\Spa{A_f[f^{-1}]}{A^{\cl}}=\R_{\X^{\ad}}\big(\frac{f}{f}\big)\mathrlap{.}
	\]
	In particular the global sections of $\mathcal{O}_W$ are the completion of $A_f[f^{-1}]$ (considered as a Huber ring), which is
	\[
	A_f[f^{-1}]^{\wedge}=\comp{A}_f\otimes_{A}A[f^{-1}]=\comp{A}_f[f^{-1}]\mathrlap{.}
	\]
	Recall that the completion of a Huber ring is given by completing its ring of definition and then tensoring with the whole ring.
	Now consider the following diagram consisting of the square overhead and the appropriate tensor and global section functors:
	\[
	\begin{tikzcd}
		&&\fMod[A]\arrow[rrr,"-\otimes_AA{[f^{-1}]}"]\arrow[ddd,swap,"-\otimes_{A}\comp{A}"]&&&\fMod[A{[f^{-1}]}]\arrow[ddd,"-\otimes_{A[f^{-1}]}\comp{A}{[f^{-1}]}"]\\
		&&&\Coh[X^{\ad}]\arrow[r]\arrow[d]\arrow[ul,"\Gamma"]\arrow[ul,"\cong",swap]&\Coh[\U]\arrow[d]\arrow[ur,"\Gamma",swap]\arrow[ur,"\cong"]\\
		&&&\Coh[\X^{\ad}]\arrow[r]\arrow[dl,"\Gamma",swap]\arrow[dl,"\cong"]&\Coh[W]\arrow[dr,"\Gamma"]\arrow[dr,"\cong",swap]\\
		&&\fMod[\comp{A}]\arrow[rrr,"-\otimes_{\comp{A}}\comp{A}{[f^{-1}]}"]&&&\fMod[\comp{A}{[f^{-1}]}]
	\end{tikzcd}
	\]
	Since $\U$ and $W$ are rational domains they are in particular affinoid, so all outer squares commute by Lemma~\ref{pullbackaffinoid}.
	The global section functors are all equivalences by Proposition~\ref{globalsec}.
	Applying Corollary~\ref{keylemmacoh} to $Z\hookrightarrow X$ yields a cartesian square
	\[
	\begin{tikzcd}
		\Coh[X]\arrow[r]\arrow[d]&\Coh[U]\arrow[d]\\
		\Coh[X^\prime]\arrow[r]&\Coh[X^\prime\times_XU]\mathrlap{,}
	\end{tikzcd}
	\]
	where we again denote $\Spec{\comp{A}}$ by $X^\prime$.
	All of these schemes are affine, so the global section functors are equivalences in this case as well.
	After composing with them this cartesian square becomes the outer square in the large diagram.
	Thus the inner square
	\[
	\begin{tikzcd}
		\Coh[X^{\ad}]\arrow[r]\arrow[d]&\Coh[\U]\arrow[d]\\
		\Coh[\X^{\ad}]\arrow[r]&\Coh[W]
	\end{tikzcd}
	\]
	in the large diagram is also cartesian, which is what we wanted to show.
	
	Now let us reduce the case of $Z$ being locally principal to the case of $X$ affine and $Z$ principal.
	Let $\{V_i\}_i$ be an open affine cover of $X$ such that $Z\cap V_i$ is the vanishing of a principal ideal in $V_i$. Let $V\subseteq X$ be some open subscheme.
	For every $i$ let $\{V^\prime_{ij}\}_{j}$ be an open affine cover of $V_i\cap V$.
	Then $Z\cap V^\prime_{ij}$ is also the vanishing of a principal ideal in $V^\prime_{ij}$, so $\{V^\prime_{ij}\}_{ij}$ is an open cover of $V$ such that the theorem holds for every $(V^\prime_{ij},Z\cap V^\prime_{ij})$ by the discussion of the affine case above.
	By Lemma~\ref{affred} the theorem holds for $(X,Z)$.
\end{proof}
\subsection{The general case}
The next step is to reduce the general case to the case of a locally principal $Z$ by generalizing Corollary~\ref{keylemmacoh} to the formal situation.

In this section let $X$ be a Noetherian affine scheme and let $Z$ be a subscheme cut out by some ideal $I$.
Let $\pi\colon\bX\to X$ be a morphism of locally Noetherian schemes and let $\widetilde{\X}$ denote the completion of $\bX$ at $\pi^{-1}(Z)=Z\times_X\bX$.

By Lemma~\ref{formalpullback} the formal scheme $\widetilde{\X}$ is the pullback $\X\times_X\bX$, so we obtain a cartesian commutative square of formal schemes
\[
\begin{tikzcd}
	\widetilde{\X}\arrow[r]\arrow[d]&\X\arrow[d]\\
	\bX\arrow[r,"\pi"]&X\mathrlap{.}
\end{tikzcd}
\]
The following lemma is the wanted generalization of Corollary~\ref{keylemmacoh}.
\begin{lemma}\label{keylemmaformal}
	Assume that $\pi$ is proper and an isomorphism over $U$.
	Then the square
	\[
	\begin{tikzcd}
		\Coh[X]\arrow[r]\arrow[d]&\Coh[\bX]\arrow[d]\\
		\Coh[\X]\arrow[r]&\Coh[\widetilde{\X}]\mathrlap{,}
	\end{tikzcd}
	\]
	which is induced by pullback along morphisms in the square above is cartesian.
\end{lemma}
\begin{remark}
	The square in the lemma is the left square of Corollary~\ref{blowupsetup} for $(X,Z)$, so this lemma will enable us to reduce to the locally principal case upon considering $\pi\colon\operatorname{Bl}_ZX\to X$.
\end{remark}
The morphism $\X\to X$ factors through $\Spec{\comp{A}}$ by the universal property of the spectrum, so we have a commutative diagram
\[
\begin{tikzcd}
	\widetilde{\X}\arrow[d,"\kappa^\prime",swap,dashed]\arrow[r]&\X\arrow[d,"\kappa"]\\
	\bX\times_X\Spec{\comp{A}}\arrow[r,"\pi^{\prime}"]\arrow[d,"i^\prime",swap]&\Spec{\comp{A}}\arrow[d,"i"]\\
	\bX\arrow[r,"\pi"]&X\mathrlap{,}
\end{tikzcd}
\]
where the outer square is the square we started with.
Let us denote $\Spec{\comp{A}}$ by $\Xs$ and $\bX\times_X\Spec{\comp{A}}$ by $\Es$.
Pulling back sheaves along maps in this diagram yields a commutative diagram
\[
\begin{tikzcd}
	\Coh[X]\arrow[r,"\pi^*"]\arrow[d,"i^*",swap]&\Coh[\bX]\arrow[d,"i^{\prime*}"]\\
	\Coh[\Xs]\arrow[r,"\pi^{\prime*}"]\arrow[d,"\kappa^*",swap]&\Coh[\Es]\arrow[d,"\kappa^{\prime*}"]\\
	\Coh[\X]\arrow[r]&\Coh[\widetilde{\X}]
\end{tikzcd}
\]
in which the large outer square is the square of Lemma~\ref{keylemmaformal}. Moreover Corollary~\ref{keylemmacoh} shows that the top square is cartesian.

Thus to show the lemma it is enough to show that $\kappa^*$ and $\kappa^{\prime*}$ are equivalences.
For this we use the following slightly technical lemma. 
\begin{lemma}\label{coherator}
	Let $B$ be a $J$-adically topologized complete Noetherian ring (with $J$ some ideal in $B$).
	Let $f\colon Y\to \Spec{B}$ be a proper map of schemes; denote the completion of $Y$ at $Y\times_{\Spec{B}}\Spec{B/J}$ by $\Y$, with structure map $\kappa\colon\Y\to Y$.
	Denote the coherator on $Y$ by
	\[
	Q\colon \Mod[\mathcal{O}_Y]\to\QCoh[Y]\mathrlap{.}
	\]
	Then
	\[
	\kappa^*\colon \Coh[Y]\to\Coh[\Y]
	\]
	is an equivalence with inverse $Q\kappa_*$.
\end{lemma}
\begin{proof}
	This is \cite[Proposition 3.1.1]{duality}.
\end{proof}
To apply this lemma to $f=\pi$ we need to know that $\widetilde{\X}$ is the completion of $\Es$ at $Z\times_X\bX$, which we check in the next lemma.
\begin{lemma}\label{compid}
	In this situation above $\widetilde{\X}$ is the completion of $\Es$ at $\pi^{\prime-1}(Z)=Z\times_X\bX$.
\end{lemma}
\begin{proof}
	We once again consider the diagram from the setup above:
	\[
	\begin{tikzcd}
		\widetilde{\X}\arrow[d,"\kappa^\prime",swap]\arrow[r]&\X\arrow[d,"\kappa"]\\
		\Es\arrow[r,"\pi^{\prime}"]\arrow[d,"i^\prime",swap]&\Xs\arrow[d,"i"]\\
		\bX\arrow[r,"\pi"]&X
	\end{tikzcd}
	\]
	The lower square is cartesian in the category of schemes, so it is also cartesian in the category of formal schemes.
	The outer square is cartesian by Lemma~\ref{formalpullback}, so the upper square is cartesian as well.
	Since $\X$ is the completion of $\Xs$ at $Z$ we can use Lemma~\ref{formalpullback} once again to see that $\widetilde{\X}$ is also the completion of $\Es$ at $\pi^{\prime-1}(Z)$.
\end{proof}
Now we are ready to prove Lemma~\ref{keylemmaformal}.
\begin{proof}[Proof of Lemma~\ref{keylemmaformal}]
	As noted earlier the outer square of the commutative diagram
	\[
	\begin{tikzcd}
		\Coh[X]\arrow[r,"\pi^*"]\arrow[d,"i^*",swap]&\Coh[\bX]\arrow[d,"i^{\prime*}"]\\
		\Coh[\Xs]\arrow[r,"\pi^{\prime*}"]\arrow[d,"\kappa^*",swap]&\Coh[\Es]\arrow[d,"\kappa^{\prime*}"]\\
		\Coh[\X]\arrow[r]&\Coh[\widetilde{\X}]
	\end{tikzcd}
	\]
	is the square we need to show to be cartesian.
	Since $\pi$ is qcqs and an isomorphism over $U$ the upper square is cartesian by Corollary~\ref{keylemmacoh}.
	Lemma~\ref{coherator} applied to $f=\id_X$ yields that $\kappa^*$ is an equivalence.
	Since $\widetilde{\X}$ is the completion of $\Es$ at $Z\times_X\bX$ by Lemma~\ref{compid} and $\pi$ is proper we can also apply Lemma~\ref{coherator} to $f=\pi$, showing that $\kappa^{\prime*}$ is an equivalence.
	Thus the outer square is cartesian, showing the claim.
\end{proof}
Putting together the reduction step and the locally principal case we have assembled everything we need to prove the main theorem.
\begin{proof}[Proof of Theorem~\ref{mainthmprop}]
	Let us first show the theorem for $(X,Z)$ where $X$ is affine. 
	Let $\pi\colon\operatorname{Bl}_ZX\to X$ be the blowup of $X$ in $Z$, denote the exceptional divisor $Z\times_X\operatorname{Bl}_ZX$ by $E$, and the completion of $\operatorname{Bl}_ZX$ at $E$ by $\E$.
	The blowup $\pi$ is proper and an isomorphism over $U$, so we can consider the commutative diagram of Corollary~\ref{blowupsetup}:
	\[
	\begin{tikzcd}
		\Coh[X]\arrow[r]\arrow[d,swap]&\Coh[\operatorname{Bl}_ZX]\arrow[r]\arrow[d]&\Coh[U]\arrow[d]\\
		\Coh[\X]\arrow[r]&\Coh[\E]\arrow[r]&\Coh[W]
	\end{tikzcd}
	\]
	Here the outer square is the square of the theorem for $(X,Z)$ and the right square is the square of the theorem for $(\operatorname{Bl}_ZX,E)$.
	
	Since $E$ is a locally principal subscheme in $\operatorname{Bl}_ZX$ the theorem holds for $(\operatorname{Bl}_ZX,E)$ by Lemma~\ref{principal}, so the right square is cartesian.
	Since $\pi$ is proper and an isomorphism over $U$ we can apply Lemma~\ref{keylemmaformal}, showing that the left square is cartesian as well.
	Thus the outer square is cartesian, which shows that the theorem holds for $(X,Z)$.
	
	Now let $(X,Z)$ be a pair with $X$ any locally Noetherian scheme.
	Let $\{V_i\}_i$ be an open affine cover of $X$ and let $V\subseteq X$ be some open subscheme.
	For every $i$ let $\{V^\prime_{ij}\}_{j}$ be an open affine cover of $V_i\cap V$.
	Then $\{V^\prime_{ij}\}_{j}$ is an open cover of $V$ such that the theorem holds for every $(V^\prime_{ij},Z\cap V^\prime_{ij})$ by the first part of the proof.
	By Lemma~\ref{affred} the theorem holds for $(X,Z)$ as well.
\end{proof}
\newpage
\section{Further discussion and outlook}
After proving a descent theorem such as Theorem~\ref{mainthmprop} one can ask which kind of properties are preserved by descending.
Unsurprisingly we obtain a similar result as we did in Corollary~\ref{keylemmafaithful} after the key lemma.
\begin{corollary}\label{preserved}
	Let $X$ be a locally Noetherian scheme, let $Z$ be a closed subscheme.
	Let $\mathcal{P}$ be a property of coherent sheaves on formal schemes such that
	\begin{itemize}
		\item 
		for every fpqc cover $\{f_i\colon Y_i\to Z\}$ of formal schemes a coherent sheaf $F$ on $Z$ has $\mathcal{P}$ if every $f_i^*F$ has $\mathcal{P}$,
		\item 
		the property $\mathcal{P}$ is stable under pullback.
	\end{itemize}
	Let $\Cohp[-]$ denote the full subcategory of $\Coh[-]$ consisting of coherent sheaves that satisfy $\mathcal{P}$.
	Then pullback yields an equivalence
	\[
	\Cohp[X]\to\Cohp[\X]\times_{\Coh[W]}\Cohp[U].
	\]
\end{corollary}
\begin{proof}
	This proof is analogous to the proof of Corollary~\ref{keylemmafaithful} when considering the square
	\[
	\begin{tikzcd}
		\Cohp[X]\arrow[r]\arrow[d,hook]&\Cp_X\arrow[d,hook]\\
		\Coh[X]\arrow[r,"\cong"]&\C_X
	\end{tikzcd}
	\]
	given by the inclusions from $\Cohp[-]$ to $\Coh[-]$.
\end{proof}
\begin{remark}
	One property that fulfills the assumptions of Corollary~\ref{preserved} is the property of being locally free, or being locally free of rank $n$.
	In particular we obtain an equivalence
	\[
	\Pic[X]\to\Pic[\X]\times_{\Coh[W]}\Pic[U]\mathrlap{.}
	\]
\end{remark}
To be able to apply the theorem one needs to understand $W$, so let us calculate $W$ for a pair $(X,Z)$ with $X$ affine.
\begin{proposition}\label{calcW}
	Let $X=\Spec{A}$ be a Noetherian affine scheme and let $Z=V(I)$ with $I=(f_1,\dots,f_n)$.
	Denote $D(f_i)$ by $U_i$; denote the preimage of $U_i$ under $s\colon \X^{\ad}\to X$ by $W_i$.
	For a multiindex $\alpha$ let $f^{\alpha}\coloneqq\prod f_j^{\alpha_j}$.
	Denote the rational subspace
	\[
	\R_{\X^{\ad}}\bigg(\frac{\{f_i\}\cup\{f^\alpha\}_{|\alpha|=k, \alpha_i=0}}{f_i}\bigg)
	\]
	by $W_{i,k}$.
	Then $W_i$ is the union over the $W_{i,k}$, so in particular
	\[
	W=\bigcup_{i\leq n}W_i=\bigcup_{i\leq n}\bigcup_{k}W_{i,k}\mathrlap{.}
	\]
\end{proposition}
\begin{proof}
	As $(\{f_i\}\cup\{f^\alpha\}_{|\alpha|=k, \alpha_i=0})$ contains $(\{f^\alpha\}_{|\alpha|=k})=I^k$ it is open in $A_I$, so the expression $\R_{\X^{\ad}}\big(\frac{\{f_i\}\cup\{f^\alpha\}_{|\alpha|=k, \alpha_i=0}}{f_i}\big)$ makes sense.
	Since $U$ is covered by the $U_i$ the preimage $W$ is covered by the $W_i$.
	
	Both $W_i$ and $\bigcup_k W_{i,k}$ are open subspaces of $\X^{\ad}$ so it is enough to check that their underlying sets agree.
	The space $W_{i,k}$ is the subspace of continuous valuations of $A_I$ that fulfill $\nu(f^\alpha)<\nu(f_i)$ for $\alpha$ with $|\alpha|=k$ and $\alpha_i=0$ as well as $\nu(f_i)\neq 0$.
	
	On the other hand the space $W_i$ is the subspace of continuous valuations of $A_I$ that do not send $f_i$ to $0$.
	Thus $\bigcup_{k}W_{i,k}$ is clearly contained in $W_i$.
	
	Let $\nu\in W_i$ be some continuous valuation. Since $\nu(f_i)\neq 0$ and $\nu$ is continuous there exists an open neighbourhood $U$ of $0$ in $A_I$ such that
	\[
	\nu(U)<\nu(f_i)\mathrlap{.}
	\]
	This neighbourhood $U$ contains some $I^k$, so in particular $\nu(f^\alpha)<\nu(f_i)$ for $|\alpha|=k$. Thus $W_i\subseteq\bigcup_{k}W_{i,k}$, and the two spaces are equal.
\end{proof}
It is moreover important to understand the intersections of the $W_i$ with each other, in the affine case we have a straightforward description.
\begin{proposition}
	With notation as in Proposition~\ref{calcW} the intersection $W_i\cap W_j$ is given by
	\[
	\bigcup_{k,l}\R_{\X^{\ad}}\bigg(\frac{\{f_if_j\}\cup\{f_if^\alpha\}\cup\{f_jf^\beta\}\cup\{f^\alpha f^\beta\}}{f_if_j}\bigg)\mathrlap{,}
	\]
	where $\alpha$ runs through indices with $|\alpha|=k$, $\alpha_i=0$, and $\beta$ runs through indices with $|\beta|=l$, $\beta_j=0$.
\end{proposition}
\begin{proof}
	This follows from the fact that for two finite sets $T_1$, $T_2$ such that $T_iA_I$ is open and $s_1$, $s_2\in A$ we have
	\[
	\R_{\X^{\ad}}\bigg(\frac{T_1}{s_1}\bigg)\cap\R_{\X^{\ad}}\bigg(\frac{T_2}{s_2}\bigg)=\R_{\X^{\ad}}\bigg(\frac{T_1T_2}{s_1s_2}\bigg)
	\]
	and the calculation of $W_i$ in the last proposition.
\end{proof}
Lastly let us say something about possible generalizations of the theorem.
One can ask whether there is a generalization involving quasicoherent modules instead of coherent modules.
Since the key lemma holds for quasicoherent modules the problem lies in the adaptation steps to sheaves on $\X$ and $W$ in chapter 3.
For sheaves on $\X$ one might be able to substitute quasicoherent modules with the subcategory of directed colimits of coherent modules $\QCohprime[\X]$ in the category of modules on $\X$, some important properties of which are developed in \cite{duality}.
By \cite[Proposition 3.1.1]{duality} the statement in Lemma~\ref{coherator} holds when replacing $\Coh[Y]$ with $\QCoh[Y]$ and $\Coh[\Y]$ with $\QCohprime[\Y]$.
Thus one also obtains an equivalence
\[
\QCoh[X]\to\QCohprime[\X]\times_{\QCohprime[\E]}\QCoh[\bX]
\]
as in Lemma~\ref{keylemmaformal}.

One might be able to also consider directed colimits of coherent modules on $W$ since the pullback functors preserve colimits, but whether the remainder of the proof can be adapted to this setting is not quite clear.

A different direction in which one can generalize the theorem is to prove a derived version with coherent modules replaced by pseudocoherent complexes.
One should be able to show a version of the key lemma in this setup with a similar proof by deriving the pushforwards and pullbacks.
There is a variant of Lemma~\ref{coherator} for derived categories in \cite[Proposition 3.3.1]{duality}, so one can likely proceed with a similar overall proof idea in this derived case.
\newpage
\bibliographystyle{alpha}
\bibliography{lit}
\end{document}